\documentclass[12pt,reqno]{amsart}
\usepackage{amsmath,amsthm,amsfonts,amssymb,latexsym,amscd,stmaryrd,amsbsy,txfonts}

\numberwithin{equation}{section}
\allowdisplaybreaks

\theoremstyle{plain}
\newtheorem{thm}{Theorem}[section]

\newtheorem{prop}[thm]{Proposition}
\newtheorem{lem}[thm]{Lemma}

\newtheorem{rem}[thm]{Remark}
\numberwithin{equation}{section}

\begin{document}

\title{Six-dimensional Painlev\'{e} systems and their particular solutions in terms of hypergeometric functions}
\date{}
\author{Takao Suzuki}
\address{Department of Mathematics, Kinki University, 3-4-1, Kowakae, Higashi-Osaka, Osaka 577-8502, Japan}
\email{suzuki@math.kindai.ac.jp}
\maketitle

\begin{abstract}
In this article, we propose a class of six-dimensional Painlev\'{e} systems given as the monodromy preserving deformations of the Fuchsian systems.
They are expressed as polynomial Hamiltonian systems of sixth order.
We also discuss their particular solutions in terms of the hypergeometric functions defined by fourth order rigid systems.


\end{abstract}

\section{Introduction}

Recently, higher order generalizations of the sixth Painlev\'{e} equation ($P_{\rm{VI}}$) has been studied from a viewpoint of the monodromy preserving deformations of Fuchsian systems.
It is shown in \cite{Kos,O1,O2} that irreducible Fuchsian systems with a fixed number of accessary parameters can be reduced to finite types of systems by using the Katz's two operations, addition and middle convolution \cite{Kat}.
It is also shown in \cite{HF} that the isomonodromy deformation equation is invariant under the Katz's two operations.
These facts allow us to construct a classification theory of Painlev\'{e} systems given as the isomonodromy deformation equations.

The Fuchsian systems with two accessary parameters were classified by Kostov \cite{Kos}.
According to it, they are reduced to the systems with the following spectral types:
\[\begin{array}{l|lll}
	\text{4 singularities} & 11,11,11,11 \\\hline
	\text{3 singularities} & 111,111,111 & 22,1111,1111 & 33,222,111111
\end{array}\]
The system with the spectral type $\{11,11,11,11\}$ gives $P_{\rm{VI}}$ as the monodromy preserving deformation \cite{Fuc}.
The other three systems have no deformation parameters, thus we can not derive isomonodromy deformation equations from them.

In general, Fuchsian systems can be classified with the aid of algorithm proposed by Oshima \cite{O1,O2}.
The systems with four accessary parameters are reduced as follows:
\[\begin{array}{l|lll}
	\text{5 singularities} & 11,11,11,11,11 \\\hline
	\text{4 singularities} & 21,21,111,111 & 31,22,22,1111 & 22,22,22,211
\end{array}\]
In addition to them, there exist nine systems which have three singularities; we do not list here.
Sakai investigated their monodoromy preserving deformations systematically and derived the four-dimensional Painlev\'{e} systems in \cite{Sak}.

An aim of this article is to investigate the monodromy preserving deformations of the Fuchsian systems with six accessary parameters.
They are reduced as follows:
\[\begin{array}{l|lll}
	\text{6 singularities} & 11,11,11,11,11,11 \\\hline
	\text{5 singularities} & 21,21,21,21,111 & 31,31,22,22,22 \\\hline
	\text{4 singularities} & 21,111,111,111 & 22,22,211,211 & 22,22,22,1111 \\
	& 31,22,211,1111 & 31,31,1111,1111 & 33,33,33,321 \\
	& 42,33,33,222 & 51,33,222,222 & 51,33,33,111111
\end{array}\]
In addition to them, there exist 24 systems which have three singularities; we do not list here.
Among those 12 systems, the following ones have been already investigated in \cite{G}, \cite{T1}, \cite{Kaw}, \cite{FISS} respectively:
\[
	\{11,11,11,11,11,11\},\quad \{31,31,1111,1111\},\quad \{33,33,33,321\},\quad \{51,33,33,111111\}.
\]
In this article, we investigate for the other eight Fuchsian systems and derive six-dimensional Painlev\'{e} systems.

It is known that $P_{\rm{VI}}$ can be expressed as the Hamiltonian system
\[
	t(t-1)\frac{dq}{dt} = \frac{\partial H_{\rm{VI}}}{\partial p},\quad
	t(t-1)\frac{dp}{dt} = -\frac{\partial H_{\rm{VI}}}{\partial q},
\]
with
\[
	H_{\rm{VI}}(\alpha_0,\alpha_1,\alpha_3,\alpha_4;q,p;t) = q(q-1)(q-t)p\left(p-\frac{\alpha_1-1}{q-t}-\frac{\alpha_3}{q-1}-\frac{\alpha_4}{q}\right) + \alpha_2(\alpha_0+\alpha_2)q,
\]
where $\alpha_0+\alpha_1+2\alpha_2+\alpha_3+\alpha_4=1$.
Such a property holds even in higher dimensional cases.
The $2n$-dimensional Painlev\'{e} systems, which have been already derived, can be expressed as $2n$-th order Hamiltonian systems
\[
	\mathcal{H}^{\mathbf{m}}:\quad
	t_i(t_i-1)\frac{\partial q_j}{\partial t_i} = \{H_i^{\mathbf{m}},q_i\},\quad
	t_i(t_i-1)\frac{\partial p_j}{\partial t_i} = \{H_i^{\mathbf{m}},p_j\}\quad (i=1,\ldots,N;j=1,\ldots,n),
\]
with the Poisson bracket defined by
\[
	\{q_i,p_j\} = -\delta_{i,j},\quad \{p_i,p_j\} = \{q_i,q_j\} = 0,
\]
where $\delta_{i,j}$ stands for the Kronecker's delta.
Here the corresponding Fuchsian systems, whose spectral type are $\mathbf{m}$, have $N+3$ singularities and $2n$ accessary parameters.
In the case $N=1$, we denote $t_1$ and $H_1$ by $t$ and $H$ respectively.
In this article, we obtain explicit formulas of the following Hamiltonians; see Section \ref{Sec:Hamiltonian}:
\[\begin{split}
	&H^{21,21,21,21,111}_i,\quad H^{31,31,22,22,22}_i\quad (i=1,2),\\
	&H^{21,111,111,111},\quad H^{31,22,211,1111},\quad H^{22,22,211,211},\quad H^{22,22,22,1111},\quad H^{42,33,33,222},\quad H^{51,33,222,222}.
\end{split}\]

\begin{rem}
The systems $\mathcal{H}^{21,21,111,111}$ and $\mathcal{H}^{31,31,1111,1111}$ were first considered by Tsuda as similarity reductions of his UC hierarchy in \cite{T1,T2}.
Independently, they were derived from similarity reductions of the Drinfeld-Sokolov hierarchy in \cite{FS3,Su1}.
The relationship between those two origins is clarified with the aid of a Laplace transformation in \cite{FS4}.
\end{rem}

\begin{rem}
The systems $\mathcal{H}^{31,22,22,1111}$ and $\mathcal{H}^{51,33,33,111111}$ were first discovered by Sasano via a generalization of the affine Weyl group symmetry and the Okamoto's initial value space for $P_{\rm{VI}}$ in \cite{Sas}.
Afterward, they were derived from similarity reductions of the Drinfeld-Sokolov hierarchy in \cite{FS2}.
\end{rem}

\begin{rem}
The system $\mathcal{H}^{51,33,33,222}$ has been already derived from a similarity reduction of the Drinfeld-Sokolov hierarchy in \cite{FS1}.
\end{rem}

The other aim of this article is to give particular solutions of the six-dimensional Painlev\'{e} systems in terms of the hypergeometric functions defined by fourth order rigid systems.
It is known that $P_{\rm{VI}}$ and the four-dimensional Painlev\'{e} systems admit particular solutions as follows:
\[\begin{array}{|l|l|l|l|}\hline
	\text{Painlev\'{e} system} & \text{Rigid system} & \text{HGF} & \text{Ref.} \\\hline
	\mathcal{H}^{11,11,11,11} & 11,11,11 & {}_2F_1 & \text{\cite{Fuc}} \\
	\mathcal{H}^{11,11,11,11,11} \ (\mathcal{H}^{21,21,21,21,21}) & 21,21,21,21 & P_3 \ (F_1) & \text{\cite{G,IKSY}} \\
	\mathcal{H}^{21,21,111,111} & 21,111,111 & {}_3F_2 & \text{\cite{Su2,T3}} \\\hline
\end{array}\]
Note that the system $\mathcal{H}^{21,21,21,21,21}$ is equivalent to $\mathcal{H}^{11,11,11,11,11}$ because the corresponding Fuchsian systems are mutually transformed by the Katz's two operations.
And, for the six-dimensional ones, we obtain the following results:
\[\begin{array}{|l|l|l|l|}\hline
	\text{Painlev\'{e} system} & \text{Rigid system} & \text{HGF} & \text{Ref.} \\\hline
	\mathcal{H}^{11,11,11,11,11,11} \ (\mathcal{H}^{31,31,31,31,31,31}) & 31,31,31,31,31 & P_4 \ (F_D) & \text{\cite{G,IKSY}} \\
	\mathcal{H}^{21,21,21,21,111} \ (\mathcal{H}^{31,31,31,22,211}) & 31,31,22,211 & \rm{II}^*_2 & \text{Sec. 5.1} \\
	\mathcal{H}^{31,31,22,22,22} & 31,22,22,22 & P_{4,4} & \text{Sec. 5.2} \\
	\mathcal{H}^{21,111,111,111} \ (\mathcal{H}^{31,211,211,211}) & 211,211,211 & \rm{II}_2 & \text{Sec. 5.3} \\
	\mathcal{H}^{31,22,211,1111} & 22,211,1111 & EO_4 & \text{Sec. 5.4} \\
	\mathcal{H}^{31,31,1111,1111} & 31,1111,1111 & {}_4F_3 & \text{\cite{Su2,T3}} \\\hline
\end{array}\]
The symbols $P_3$, $P_4$ and $P_{4,4}$ stand for the Jordan-Pochhammer family and its generalization (cf. \cite{O2}).
And $EO_4$ stands for the even-four hypergeometric function which is in Simpson's list \cite{Sim}.
Moreover $\rm{II}_2$ and $\rm{II}^*_2$ are in Yokoyama's list \cite{Y}.

\begin{rem}
As is seen in above, for $n=1,2,3$, the $2n$-dimensional Painlev\'{e} system $\mathcal{H}^{n1,m_2,\ldots,m_{N+3}}$ admits a particular solution in terms of the $(n+1)$-th order rigid system with a spectral type $\{m_2,\ldots,m_{N+3}\}$.
Such a relationship is satisfied for a more general case{\rm;} see Section \ref{Sec:Solutions}.
\end{rem}

This article is organized as follows.
In Section \ref{Sec:Hamiltonian} and Appendix \ref{App:Hamiltonian_6dim}, we give explicit formulas of the six-dimensional Painlev\'{e} systems.
In Section \ref{Sec:Schlesinger}, we recall the Schlesinger system and its Poisson structure.
In Section \ref{Sec:Derivation}, we discuss derivations of the Painlev\'{e} systems from the Schlesinger systems.
In Section \ref{Sec:Solutions}, we give particular solutions of the six-dimensional Painlev\'{e} systems in terms of the hypergeometric functions.
In Appendix \ref{App:Hamiltonian_4dim}, we recall the four-dimensional Painlev\'{e} systems which have been classified by Sakai.

\section{List of Hamiltonians obtained in this article}\label{Sec:Hamiltonian}

In this section, we give explicit formulas of the following Hamiltonians:
\[\begin{split}
	&H^{21,21,21,21,111}_i,\quad H^{31,31,22,22,22}_i\quad (i=1,2),\\
	&H^{21,111,111,111},\quad H^{31,22,211,1111},\quad H^{22,22,211,211},\quad H^{22,22,22,1111},\quad H^{42,33,33,222},\quad H^{51,33,222,222}.
\end{split}\]
The following ones are given in Appendix \ref{App:Hamiltonian_6dim}:
\[
	H^{11,11,11,11,11,11}_i,\quad (i=1,2,3),\quad H^{31,31,1111,1111},\quad H^{33,33,33,321},\quad H^{51,33,33,111111}.
\]
And the following ones are given in Appendix \ref{App:Hamiltonian_4dim}:
\[
	H^{11,11,11,11,11}_i,\quad (i=1,2),\quad H^{21,21,111,111},\quad H^{22,22,22,211},\quad H^{31,22,22,1111}.
\]

\subsection{Spectral type $21,21,21,21,111$}

\[\begin{split}
	&H^{21,21,21,21,111}_1\\
	&= H^{11,11,11,11,11}_1 + \frac{t_1-1}{t_1-t_2}q_3(q_3-t_2)(q_3-t_1)p_3\left(p_3-\frac{\alpha_0+\alpha_3+\alpha_5}{q_3-t_1}-\frac{\alpha_0+\alpha_4+\rho_3}{q_3-t_2}-\frac{\alpha_1}{q_3}\right)\\
	&\quad - (\alpha_0+\alpha_5+1)q_3p_3 - \frac{t_1-1}{t_1-t_2}\alpha_2\rho_3q_3 + \frac{1}{t_2}q_2q_3(q_1p_1+q_2p_2+\alpha_0+\alpha_5+1)(q_3p_3-\rho_3)\\
	&\quad + q_1q_3p_3(q_1p_1+q_2p_2+\alpha_0+\alpha_5+1) + \frac{1}{t_2(t_1-t_2)}q_2q_3\{t_1(t_2-1)p_1-(t_1-1)t_2p_2\}(q_3p_3-\rho_3)\\
	&\quad + \frac{1}{t_1-t_2}q_3p_3\{-(t_1^2+t_1-2t_2)q_1p_1-t_1(t_1-1)q_2p_1+(t_1-1)t_2q_1p_2+t_1(t_2-1)q_2p_2\}\\
	&\quad + \frac{t_1}{t_1-t_2}p_3\{(t_1-1)t_2q_1p_1-(t_1-1)t_2q_1p_2+(t_1-t_2)q_3p_1\},
\end{split}\]
and
\[\begin{split}
	&H^{21,21,21,21,111}_2\\
	&= H^{11,11,11,11,11}_2 + \rho_3q_2\{(q_2-1)p_2+q_1p_1+\alpha_0+\alpha_5+1\}\\
	&\quad + \frac{t_2-1}{t_2-t_1}q_3(q_3-t_1)(q_3-t_2)p_3\left(p_3-\frac{\alpha_0+\alpha_4+\rho_3-1}{q_3-t_2}-\frac{\alpha_0+\alpha_3+\alpha_5+1}{q_3-t_1}-\frac{\alpha_1}{q_3}\right)\\
	&\quad - (\alpha_0+\alpha_5+1)q_3p_3 - \frac{t_2-1}{t_2-t_1}\alpha_2\rho_3q_3 + q_2q_3p_3(q_2p_2+q_1p_1+\alpha_0+\alpha_5+1)\\
	&\quad - \frac{t_2-1}{t_2-t_1}q_2q_3(p_2-p_1)(q_3p_3-\rho_3) + t_2q_1p_3(q_2p_2+q_1p_1+\alpha_0+\alpha_5+1)\\
	&\quad + \frac{1}{t_2-t_1}q_3p_3\{(t_2^2+t_2-2t_1)q_2p_2+t_2(t_2-1)q_1p_2-(t_2-1)t_1q_2p_1-t_2(t_1-1)q_1p_1\}\\
	&\quad + \frac{t_2}{t_2-t_1}p_3\{-t_2(t_1-1)q_1p_2+(t_2-1)t_1q_1p_1-(t_2-t_1)q_3p_2\},
\end{split}\]
where
\[
	\alpha_0 = -\rho_2,\quad \alpha_1 = -\theta_1,\quad \alpha_2 = -\theta_2-\rho_3,\quad \alpha_3 = -\theta_4,\quad \alpha_4 = -\rho_1+\rho_2+1,\quad \alpha_5 = -\theta_3-1,
\]
and $\theta_1+\theta_2+\theta_3+\theta_4+\rho_1+\rho_2+\rho_3=0$.

\begin{rem}
The system $\mathcal{H}^{21,21,21,21,111}$ reduces to $\mathcal{H}^{11,11,11,11,11}$ via a specialization $p_3=\rho_3=0$.
\end{rem}

\subsection{Spectral type $31,31,22,22,22$}

\[\begin{split}
	H^{31,31,22,22,22}_i &= H^{11,11,11,11,11}_i - (\alpha_1-\alpha_2)\frac{t_1}{t_1-t_2}q_2\{(t_i-1)p_i-(t_2-1)p_{i+1}\}\\
	&\quad + \delta_{i,2}(\alpha_1-\alpha_2)t_2(q_2-1)p_2 - (t_i+1)q_3^2p_3^2 + \{\alpha_1+\alpha_3-1+(\alpha_1+\alpha_4)t_i\}q_3p_3\\
	&\quad + q_iq_3p_3(2q_1p_1+2q_2p_2+q_3p_3+2\alpha_0+\alpha_5+1)\\
	&\quad - q_3p_{i+1}(2q_ip_i+q_{i+1}p_{i+1}+2q_3p_3+2\alpha_0+\alpha_2+\alpha_5+1)\\
	&\quad - 2(t_i+1)q_iq_3p_ip_3 + t_iq_{i+1}p_3(q_3p_3-\alpha_1+\alpha_2)\\
	&\quad + 2t_iq_3p_ip_3 + \frac{1}{t_i-t_{i+1}}q_3\{t_i(t_i-1)p_i^2-2t_i(t_{i+1}-1)p_1p_2+(t_i-1)t_{i+1}p_{i+1}^2\},
\end{split}\]
for $i\in\mathbb{Z}/2\mathbb{Z}$, where
\[
	\alpha_0 = \frac{\theta_1-\theta_2-2\rho_2}{2},\quad \alpha_1 = -\theta_1,\quad \alpha_2 = \frac{-\theta_1+\theta_2}{2},\quad \alpha_3 = -\theta_4,\quad \alpha_4 = -\rho_1+\rho_2+1,\quad \alpha_5 = -\theta_3-1,
\]
and $\theta_1+\theta_2+2\theta_3+2\theta_4+2\rho_1+2\rho_2=0$.

\begin{rem}
The system $\mathcal{H}^{31,31,22,22,22}$ reduces to $\mathcal{H}^{11,11,11,11,11}$ with $\alpha_1=\alpha_2$ via a specialization $q_3=\theta_1+\theta_2=0$.
\end{rem}

\subsection{Spectral type $21,111,111,111$}

\[\begin{split}
	H^{21,111,111,111} &= H^{21,21,111,111} + q_1q_3p_3(2q_1p_1+q_2p_2+q_3p_3+\alpha_1+\eta)\\
	&\quad - q_2q_3p_3(q_1p_1+q_3p_3-\alpha_5+\eta) + tq_2q_3p_1(q_1p_1+q_2p_2+q_3p_3+\eta)\\
	&\quad + (q_1-q_2)(tp_1-p_3)(q_3p_3+\theta_{2,1}) - q_1q_3p_3(tp_1+p_2)\\
	&\quad - tq_3p_1(q_1p_1+q_2p_2-\alpha_3+\eta),
\end{split}\]
where
\[\begin{split}
	&\alpha_0 = \rho_2,\quad \alpha_1 = -\theta_{3,2}-\rho_2,\quad \alpha_2 = -\theta_{3,1}+\theta_{3,2},\quad \alpha_3 = \theta_{2,1}+\theta_{3,1}+\rho_1,\\
	&\alpha_4 = -\rho_1+\rho_3+1,\quad \alpha_5 = -\theta_{2,1}-\rho_3,\quad \eta = \theta_1+\theta_{2,1}+\theta_{3,1}+\rho_1,
\end{split}\]
and $\theta_1+\theta_{2,1}+\theta_{2,2}+\theta_{3,1}+\theta_{3,2}+\rho_1+\rho_2+\rho_3=0$.

\begin{rem}
The system $\mathcal{H}^{21,111,111,111}$ reduces to $\mathcal{H}^{21,21,111,111}$ via a specialization $q_3=\theta_{2,1}=0$.
\end{rem}

\subsection{Spectral type $31,22,211,1111$}

\[\begin{split}
	H^{31,22,211,1111} &= H^{21,21,111,111} - \rho_4(\alpha_5-\eta+\rho_4)q_2\\
	&\quad + H_{\rm{VI}}(-\alpha_5+\eta-2\rho_4,\alpha_0,\alpha_2+\alpha_3+\alpha_4+\eta,\alpha_1;q_3,p_3;t)\\
	&\quad - \frac{1}{t}q_1q_3(q_3p_3-\rho_4)(q_3p_3+\alpha_5-\eta+\rho_4)\\
	&\quad - q_1q_3p_3(2q_1p_1+q_2p_2-q_3p_3+\alpha_1+\eta) + q_2q_3p_3(2q_1p_1+q_2p_2+q_3p_3+\alpha_5)\\
	&\quad + tq_2p_3(q_1p_1-q_2p_2-q_3p_3+\alpha_1) + q_3p_3\{(t+1)q_1p_1-2tq_2p_1+q_1p_2-q_2p_2\}\\
	&\quad - tp_1p_3(q_1-tq_2),
\end{split}\]
where
\[\begin{split}
	&\alpha_0 = \theta_{3,2},\quad \alpha_1 = -\theta_2-\theta_{3,2}-\rho_2-\rho_4,\quad \alpha_2 = \theta_2+\theta_{3,1}+\rho_2+\rho_4,\\
	&\alpha_3 = -\theta_2-\theta_{3,1}-\rho_1-\rho_4,\quad \alpha_4 = \rho_1-\rho_3+1,\quad \alpha_5 = \theta_2+\rho_3+\rho_4,\quad \eta = \rho_3-\rho_4,
\end{split}\]
and $\theta_1+2\theta_2+\theta_{3,1}+\theta_{3,2}+\rho_1+\rho_2+\rho_3+\rho_4=0$.

\begin{rem}
The system $\mathcal{H}^{31,22,211,1111}$ reduces to $\mathcal{H}^{21,21,111,111}$ via a specialization $p_3=\rho_4=0$.
It also reduces to $\mathcal{H}^{31,22,22,1111}$ via a specialization $(q_1-t)p_1-q_3p_3-\alpha_2=\alpha_1+\alpha_2=0$.
\end{rem}

\subsection{Spectral type $22,22,211,211$}

\[\begin{split}
	H^{22,22,211,211} &= \mathrm{tr}\{Q(Q-1)P(Q-t)P-(\alpha_1-1)Q(Q-1)P-\alpha_3Q(Q-t)P\\
	&\quad -\alpha_4(Q-1)(Q-t)P+\alpha_2(\alpha_0+\alpha_2)Q+(t-1)R(QP-\alpha_4)\},
\end{split}\]
with
\[
	P = -\frac{1}{t}\begin{pmatrix}p_1&p_2\\q_2p_2-\alpha_2-\alpha_5&p_3\end{pmatrix},\quad
	Q = -t\begin{pmatrix}q_1&1\\q_2&q_3\end{pmatrix},\quad
	R = \begin{pmatrix}0&-(q_1-q_3)p_2+p_1-p_3\\0&\alpha_1+\theta_{3,2}-1\end{pmatrix},
\]
where
\[\begin{split}
	&\alpha_0 = -\theta_2,\quad \alpha_1 = -\theta_{3,1}+1,\quad \alpha_2 = -\rho_3,\\
	&\alpha_3 = \theta_1+\theta_2+\theta_{3,1}+2\rho_3,\quad \alpha_4 = -\theta_1,\quad \alpha_5 = -\theta_1-\theta_2-\theta_{3,2}-\rho_2,
\end{split}\]
and $2\theta_1+2\theta_2+\theta_{3,1}+\theta_{3,2}+\rho_1+\rho_2+2\rho_3=0$.

\begin{rem}
The system $\mathcal{H}^{22,22,211,211}$ reduces to $\mathcal{H}^{22,22,22,211}$ via a specialization $R=O$.
\end{rem}

\subsection{Spectral type $22,22,22,1111$}

\[\begin{split}
	H^{22,22,22,1111} &= \mathrm{tr}\{Q(Q-1)P(Q-t)P-(\alpha_1-1)Q(Q-1)P-\alpha_3Q(Q-t)P\\
	&\quad -\alpha_4(Q-1)(Q-t)P+\alpha_2(\alpha_0+\alpha_2)Q-R(QP+PQ)(Q-t)-(\alpha_0+\alpha_2-\rho_4)RQ\},
\end{split}\]
with
\[
	P = -\frac{1}{t}\begin{pmatrix}p_1&p_2\\q_2p_2-\alpha_2-\alpha_5&p_3\end{pmatrix},\quad
	Q = -t\begin{pmatrix}q_3&1\\q_2&q_1\end{pmatrix},\quad
	R = \begin{pmatrix}\alpha_2+\rho_4&-(q_1-q_3)p_2+p_1-p_3\\0&0\end{pmatrix},
\]
where
\[\begin{split}
	&\alpha_0 = -\theta_2,\quad \alpha_1 = -\theta_3+1,\quad \alpha_2 = -\rho_3,\\
	&\alpha_3 = \theta_1+\theta_2+\theta_3+2\rho_3,\quad \alpha_4 = -\theta_1,\quad \alpha_5 = -\theta_1-\theta_2-\theta_3-\rho_2,
\end{split}\]
and $2\theta_1+2\theta_2+2\theta_3+\rho_1+\rho_2+\rho_3+\rho_4=0$.

\begin{rem}
The system $\mathcal{H}^{22,22,22,1111}$ reduces to $\mathcal{H}^{22,22,22,211}$ via a specialization $R=O$.
\end{rem}

\subsection{Spectral type $42,33,33,222$}

\[\begin{split}
	&H^{42,33,33,222}\\
	&= H_{\rm{VI}}(\theta_1+\theta_2+\rho_2,\theta_1+2\theta_2+\theta_3+\rho_2+2\rho_3+1,-\theta_1-\theta_3,-\theta_1-\theta_2-2\rho_2-2\rho_3;q_1,p_1;t)\\
	&\quad + q_1p_1(tp_1+\rho_2) - \frac{\theta_1-\theta_2+2\rho_2}{2}q_2p_2\\
	&\quad + H_{\rm{VI}}(-\frac{\theta_1-\theta_2+2\rho_2}{2},-\frac{\theta_1+\theta_2+2\theta_3+2\rho_2}{2},\theta_1+2\theta_2+\theta_3+2\rho_3+1,2\rho_2;q_2,p_2;t)\\
	&\quad - 2(t+1)q_3^2p_3^2 + \{-\theta_1-\theta_2-2\theta_3-1+(3\theta_1+3\theta_2+2\theta_3+4\rho_2+4\rho_3+1)t\}q_3p_3\\
	&\quad - q_1q_3(q_1p_1-\theta_2)(q_1p_1+\theta_1+\rho_2)\\
	&\quad - q_1q_2p_1(q_1p_1+\frac{\theta_1-\theta_2+2\rho_2}{2}) + q_2q_3p_3(2q_2p_2+q_3p_3-\frac{\theta_1+3\theta_2+2\rho_2+4\rho_3}{2})\\
	&\quad + tq_3p_2(-2q_1p_1+3q_2p_2+4q_3p_3-2\theta_1-2\theta_2-2\rho_2-4\rho_3)\\
	&\quad + q_2p_3(2q_1p_1-q_3p_3+\frac{\theta_1-\theta_2+2\rho_2}{2}) + 2(t+1)q_3p_3(q_1p_1-q_2p_2)\\
	&\quad - 2t\{(t+1)q_3p_2^2+p_1p_2(q_1+q_2)+q_3p_3(p_1-2p_2)\} + 2tp_1\{(t+1)p_2-p_3\},
\end{split}\]
where $3\theta_1+3\theta_2+2\theta_3+2\rho_1+2\rho_2+2\rho_3=0$.

\subsection{Spectral type $51,33,222,222$}

\[\begin{split}
	H^{51,33,222,222} &= H_{\rm{VI}}(\alpha_3,-\alpha_1-2\alpha_2-2\alpha_3+1,\alpha_1,\alpha_3;q_1,p_1;t)\\
	&\quad + H_{\rm{VI}}(\alpha_3,-2\alpha_3-2\alpha_4-\alpha_5+1,\alpha_5,\alpha_3;q_2,p_2;t)\\
	&\quad + H_{\rm{VI}}(\alpha_3,-\alpha_0-2\alpha_3-2\alpha_6+1,\alpha_0,\alpha_3;q_3,p_3;t)\\
	&\quad + \{(q_1-1)p_1+\alpha_2\}\{(q_2-1)p_2+\alpha_4\}(q_1q_2+t)\\
	&\quad + \{(q_2-1)p_2+\alpha_4\}\{(q_3-1)p_3+\alpha_6\}(q_2q_3+t)\\
	&\quad + \{(q_3-1)p_3+\alpha_6\}\{(q_1-1)p_1+\alpha_2\}(q_3q_1+t),
\end{split}\]
where 
\[\begin{split}
	&\alpha_0 = \theta_1+\theta_{3,1}+\theta_{3,2}+\rho_2+\rho_3+1,\quad
	\alpha_1 = \theta_1+\rho_2,\quad
	\alpha_2 = -\theta_1,\quad
	\alpha_3 = -\rho_3,\\
	&\alpha_4 = -\frac{\theta_1}{2}-\frac{\theta_2}{2}-\theta_{3,2}-\rho_2,\quad
	\alpha_5 = -\theta_{3,1}+\theta_{3,2},\quad
	\alpha_6 = \frac{\theta_1}{2}+\frac{\theta_2}{2}+\rho_3,
\end{split}\]
and $3\theta_1+\theta_2+2\theta_{3,1}+2\theta_{3,2}+2\rho_1+2\rho_2+2\rho_3=0$.

\section{Schlesinger system}\label{Sec:Schlesinger}

In this section, we recall the Schlesinger system and its Poisson structure following the previous work \cite{JMMS,JMU,Sak,Sch}.

Let $A_1,\ldots,A_{N+2}\in M_{L}(\mathbb{C})$.
We consider a Fuchsian system on $\mathbb{P}^1(\mathbb{C})$
\begin{equation}\label{Eq:Fuchs}
	\frac{\partial}{\partial x}Y(x) = \sum_{i=1}^{N+2}\frac{A_i}{x-t_i}Y(x),
\end{equation}
with regular singularities $x=t_1,\ldots,t_N,t_{N+1}=1,t_{N+2}=0,\infty$.
Here we assume that each $A_i$ $(i=1,\ldots,N+2)$ can be diagonalized and $A_{\infty}:=-\sum_{i=1}^{N+2}A_i$ is a diagonal matrix.
Then the monodromy preserving deformation of \eqref{Eq:Fuchs} gives a Schlesinger system
\[
	\frac{\partial A_j}{\partial t_i} = \frac{[A_i,A_j]}{t_i-t_j},\quad
	\frac{\partial A_i}{\partial t_i} = -\sum_{j=1; j\neq i}^{N+2}\frac{[A_i,A_j]}{t_i-t_j}\quad (i=1,\ldots,N; j=1,\ldots,N+2; j\neq i).
\]
It is also expressed as a Hamiltonian system
\begin{equation}\label{Eq:Schlesinger}
	\frac{\partial A_j}{\partial t_i} = \{H_i,A_j\},\quad
	H_i = \sum_{j=1; j\neq i}^{N+2}\frac{\mathrm{tr}A_iA_j}{t_i-t_j}\quad (i=1,\ldots,N; j=1,\ldots,N+2),
\end{equation}
with a Poisson bracket defined by
\begin{equation}\label{Eq:Poisson}
	\{(A_i)_{k,l},(A_j)_{r,s}\} = \delta_{i,j}\{\delta_{r,l}(A_i)_{k,s}-\delta_{k,s}(A_i)_{r,l}\}.
\end{equation}

Thanks to the method established in \cite{JMMS}, the Schlesinger system can be rewritten to a canonical Hamiltonian system.
Consider a decomposition of matrices $A_i$ as
\[
	A_i = B_iC_i,\quad
	B_i = \left(b^{(i)}_{k,l}\right)_{k,l} \in M_{L,\mathrm{rank}A_i}(\mathbb{C}),\quad
	C_i = \left(c^{(i)}_{l,k}\right)_{l,k} \in M_{\mathrm{rank}A_i,L}(\mathbb{C}).
\]
Then the variables $b^{(i)}_{k,l},c^{(i)}_{l,k}$ can be regarded as canonical ones.
In fact, the Poisson bracket
\[
	\{b^{(i)}_{k,l},c^{(i)}_{l,k}\} = -1,\quad
	\{\text{otherwise}\} = 0,
\]
implies the above one \eqref{Eq:Poisson}.
In terms of those variables, the system \eqref{Eq:Schlesinger} is expressed as a Hamiltonian system
\begin{equation}\label{Eq:Sch_Ham}
	\frac{\partial B_j}{\partial t_i} = \{H_i,B_j\},\quad
	\frac{\partial C_j}{\partial t_i} = \{H_i,C_j\},\quad
	H_i = \sum_{j=1; j\neq i}^{N+2}\frac{\mathrm{tr}B_iC_iB_jC_j}{t_i-t_j}\quad (i=1,\ldots,N; j=1,\ldots,N+2),
\end{equation}
with a symplectic form
\[
	\omega = \sum_{i=1}^{N+2}\mathrm{tr}(dB_i\wedge dC_i) = \sum_{i=1}^{N+2}\sum_{k=1}^{L}\sum_{l=1}^{\mathrm{rank}A_i}db^{(i)}_{k,l}\wedge dc^{(i)}_{l,k}.
\]

It remains to find a canonical variables which is suitable for the number of accessory parameters of \eqref{Eq:Fuchs}.
We denote the multiplicity data of eigenvalues of $A_1,\ldots,A_{N+2},A_{\infty}$, called {\it a spectral type}, by a $(N+3)$-tuples of partitions of natural number $L$
\[
	\left\{(m_{1,1},\ldots,m_{1,l_1}),\ldots,(m_{N+2,1},\ldots,m_{N+2,l_{N+2}}),(m_{\infty,1},\ldots,m_{\infty,l_{\infty}})\right\}.
\]
Then the number of accessory parameters of \eqref{Eq:Fuchs} is given by
\[
	(N+1)L^2 - \sum_{i=1}^{N+2}\sum_{j=1}^{l_i}m_{i,j}^2 - \sum_{j=1}^{l_{\infty}}m_{\infty,j}^2 + 2.
\]
And it is generally less than the dimension of a space of matrices $(B_1,C_1,\ldots,B_{N+2},C_{N+2})$.
When we reduce the number of dependent variables of \eqref{Eq:Sch_Ham} to the suitable one, the following proposition plays an important role.
\begin{prop}[\cite{Sak}]\label{Prop:Sch_Ham_Transf}
Let $P\in GL_L(\mathbb{C})$ and $Q_j\in GL_{\mathrm{rank}A_j}(\mathbb{C})$ $(j=1,\ldots,N+2)$.
\begin{enumerate}
\item
If $\mathrm{tr}\{A_{\infty}(dP)P^{-1}\wedge(dP)P^{-1}\}=0$, then $\omega=\sum_{i=1}^{N+2}\mathrm{tr}(dP^{-1}B_i\wedge dC_iP)$.
\item
If $dQ_jC_jB_jQ_j^{-1}=0$ and $\mathrm{tr}\{Q_jC_jB_jQ_j^{-1}(dQ_j)Q_j^{-1}\wedge(dQ_j)Q_j^{-1}\}=0$, then $\omega=\mathrm{tr}(dB_jQ_j^{-1}\wedge dQ_jC_j)+\sum_{i=1; i\neq j}^{N+2}\mathrm{tr}(dB_i\wedge dC_i)$.
\end{enumerate}
\end{prop}

\section{Derivation of the Painlev\'{e} system}\label{Sec:Derivation}

In this section, we derive six-dimensional Painlev\'{e} systems from the Schlesinger system \eqref{Eq:Sch_Ham} associated with the following spectral types:
\[\begin{array}{llll}
	\{21,21,21,21,111\}, & \{31,31,22,22,22\}, & \{21,111,111,111\}, & \{31,22,211,1111\}, \\
	\{22,22,211,211\}, & \{22,22,22,1111\}, & \{42,33,33,222\}, & \{51,33,222,222\}.
\end{array}\]

\subsection{Spectral type $21,21,21,21,111$}\label{Sec:Derivation_21_21_21_21_111}

We consider a Fuchsian system
\begin{equation}\label{Eq:Fuchs_21_21_21_21_111}
	\frac{\partial}{\partial x}Y(x) = \left(\frac{A_1}{x-t_1}+\frac{A_2}{x-t_2}+\frac{A_3}{x-1}+\frac{A_4}{x}\right)Y(x),
\end{equation}
with a Riemann scheme
\[
	\left\{\begin{array}{ccccc}
		x=t_1 & x=t_2 & x=1 & x=0 & x=\infty \\
		\theta_1 & \theta_2 & \theta_3 & \theta_4 & \rho_1 \\
		0 & 0 & 0 & 0 & \rho_2 \\
		0 & 0 & 0 & 0 & \rho_3
	\end{array}\right\}.
\]
Note that a Fuchsian relation $\theta_1+\theta_2+\theta_3+\theta_4+\rho_1+\rho_2+\rho_3=0$ is satisfied.
The residue matrices are expressed as
\[
	A_i = B_iC_i,\quad
	B_i = \begin{pmatrix}b^{(i)}_1\\b^{(i)}_2\\b^{(i)}_3\end{pmatrix},\quad
	C_i = \begin{pmatrix}c^{(i)}_1&c^{(i)}_2&c^{(i)}_3\end{pmatrix}\quad (i=1,2,3,4),
\]
where $b^{(i)}_1c^{(i)}_1+b^{(i)}_2c^{(i)}_2+b^{(i)}_3c^{(i)}_3=\theta_i$, and
\[
	A_{\infty} := -\sum_{i=1}^{4}A_i = \begin{pmatrix}\rho_1&0&0\\0&\rho_2&0\\0&0&\rho_3\end{pmatrix}.
\]

Under the Schlesinger system \eqref{Eq:Sch_Ham} associated with the Fuchsian one \eqref{Eq:Fuchs_21_21_21_21_111}, we consider a gauge transformation
\[
	\widetilde{A}_i = \widetilde{B}_i\widetilde{C}_i,\quad
	\widetilde{B}_i = P^{-1}B_iQ_i^{-1},\quad
	\widetilde{C}_i = Q_iC_iP\quad (i=1,2,3,4),
\]
where
\[
	P = \begin{pmatrix}1&0&0\\\frac{b^{(4)}_2}{b^{(4)}_1}&\frac{1}{b^{(4)}_1c^{(4)}_2}&0\\\frac{b^{(4)}_3}{b^{(4)}_1}&0&\frac{1}{b^{(4)}_1c^{(4)}_3}\end{pmatrix},\quad
	Q_i = \begin{pmatrix}b^{(i)}_1\end{pmatrix}\quad (i=1,2,3,4).
\]
Then the residue matrices are transformed into
\[\begin{split}
	&\widetilde{A}_1 = \begin{pmatrix}1\\b_1\\b_2\end{pmatrix}\begin{pmatrix}\theta_1-b_1c_1-b_2c_2&c_1&c_2\end{pmatrix},\quad
	\widetilde{A}_2 = \begin{pmatrix}1\\b_3\\b_4\end{pmatrix}\begin{pmatrix}\theta_2-b_3c_3-b_4c_4&c_3&c_4\end{pmatrix},\\
	&\widetilde{A}_3 = \begin{pmatrix}1\\a_1\\a_2\end{pmatrix}\begin{pmatrix}\theta_3-a_1-a_2&1&1\end{pmatrix},\quad
	\widetilde{A}_4 = \begin{pmatrix}1\\0\\0\end{pmatrix}\begin{pmatrix}\theta_4&a_3&a_4\end{pmatrix},
\end{split}\]
and
\begin{equation}\label{Eq:Residue_21_21_21_21_111}
	\widetilde{A}_{\infty} := -\sum_{i=1}^{4}\widetilde{A}_i = \begin{pmatrix}\rho_1&0&0\\a_5&\rho_2&0\\a_6&0&\rho_3\end{pmatrix}.
\end{equation}
Note that each component is rational in $(b^{(i)}_j,c^{(i)}_j)$; we do not give its explicit formula here.
Furthermore, the relation \eqref{Eq:Residue_21_21_21_21_111} implies
\[\begin{split}
	a_1 &= -\rho_2 - b_1c_1 - b_3c_3,\\
	a_2 &= -\rho_3 - b_2c_2 - b_4c_4,\\
	a_3 &= -c_1 - c_3 - 1,\\
	a_4 &= -c_2 - c_4 - 1,\\
	a_5 &= -b_1(\theta_1-b_1c_1-b_2c_2) - b_3(\theta_2-b_3c_3-b_4c_4)\\
	&\quad + (\rho_2+b_1c_1+b_3c_3)(\theta_3+\rho_2+\rho_3+b_1c_1+b_2c_2+b_3c_3+b_4c_4),\\
	a_6 &= -b_2(\theta_1-b_1c_1-b_2c_2) - b_4(\theta_2-b_3c_3-b_4c_4)\\
	&\quad + (\rho_3+b_2c_2+b_4c_4)(\theta_3+\rho_2+\rho_3+b_1c_1+b_2c_2+b_3c_3+b_4c_4),
\end{split}\]
and
\begin{equation}\label{Eq:DepVar_21_21_21_21_111}
	b_1(c_1-c_2) + b_3(c_3-c_4) + \rho_2 = 0,\quad
	b_2(c_2-c_1) + b_4(c_4-c_3) + \rho_3 = 0.
\end{equation}
Hence the components of $\widetilde{A}_1,\ldots,\widetilde{A}_4$ turn out to be polynomials in $(b_j,c_j)$.
Then, thanks to Proposition \ref{Prop:Sch_Ham_Transf}, we obtain
\begin{prop}
The dependent variables $b_j,c_j$ $(j=1,2,3,4)$ satisfy a Hamiltonian system of eighth order
\begin{equation}\begin{split}\label{Eq:Sch_Ham_21_21_21_21_111}
	&\frac{\partial b_j}{\partial t_i} = \{H_i,b_j\},\quad
	\frac{\partial c_j}{\partial t_i} = \{H_i,c_j\}\quad (i=1,2),\\
	&H_1 = \frac{\mathrm{tr}\widetilde{A}_1\widetilde{A}_2}{t_1-t_2} + \frac{\mathrm{tr}\widetilde{A}_1\widetilde{A}_3}{t_1-1} + \frac{\mathrm{tr}\widetilde{A}_1\widetilde{A}_4}{t_1},\quad
	H_2 = \frac{\mathrm{tr}\widetilde{A}_2\widetilde{A}_1}{t_2-t_1} + \frac{\mathrm{tr}\widetilde{A}_2\widetilde{A}_3}{t_2-1} + \frac{\mathrm{tr}\widetilde{A}_2\widetilde{A}_4}{t_2},
\end{split}\end{equation}
with a symplectic form
\begin{equation}\label{Eq:Simplectic_21_21_21_21_111}
	\omega = \sum_{j=1}^{4}db_j\wedge dc_j,
\end{equation}
and the relation \eqref{Eq:DepVar_21_21_21_21_111}.
\end{prop}

We next reduce the Hamiltonian system \eqref{Eq:Sch_Ham_21_21_21_21_111} to the one of sixth order.
Substituting the second relation of \eqref{Eq:DepVar_21_21_21_21_111} to \eqref{Eq:Simplectic_21_21_21_21_111}, we obtain
\[\begin{split}
	\omega &= d(b_1+b_2)\wedge dc_1 + db_2\wedge d(c_2-c_1) + d(b_3+b_4)\wedge dc_3 + db_4\wedge d(c_4-c_3)\\
	&= d(b_1+b_2)\wedge dc_1 + db_2\wedge d(c_2-c_1) + d(b_3+b_4)\wedge dc_3 - d\frac{b_2(c_2-c_1)}{c_4-c_3}\wedge d(c_4-c_3)\\
	&= d(b_1+b_2)\wedge dc_1 + db_2(c_4-c_3)\wedge d\frac{c_2-c_1}{c_4-c_3} + d(b_3+b_4)\wedge dc_3.
\end{split}\]
Hence we can take a six-dimensional canonical coordinate system by
\[\begin{split}
	&\frac{q_1}{t_1} = -c_1,\quad
	t_1p_1 = b_1 + b_2,\quad
	\frac{q_2}{t_2} = -c_3,\quad
	t_2p_2 = b_3 + b_4,\\
	&\frac{q_3}{t_1} = -\frac{c_2-c_1}{c_4-c_3},\quad
	t_1p_3 = b_2(c_4-c_3).
\end{split}\]
Let
\[
	\widetilde{H}_1 = H_1 + \frac{q_1p_1}{t_1} + \frac{q_3p_3}{t_1},\quad
	\widetilde{H}_2 = H_2 + \frac{q_2p_2}{t_2}.
\]
Then it is easy to verify that the Hamiltonian $\widetilde{H}_i$ is just equivalent to the one $H^{21,21,21,21,111}_i$, which was given in Section \ref{Sec:Hamiltonian}, for each $i=1,2$.
Note that the variables $b_j,c_j$ $(j=1,\ldots,4)$ are described in terms of the canonical coordinates as
\[\begin{split}
	&b_1 = t_1p_1 - t_1p_3\frac{t_2p_2-q_3p_1}{\rho_2-\rho_3},\quad
	c_1 = -\frac{q_1}{t_1},\\
	&b_2 = t_1p_3\frac{t_2p_2-q_3p_1}{\rho_2-\rho_3},\quad
	c_2 = -\frac{q_3}{t_1}\frac{\rho_2-\rho_3}{t_2p_2-p_1q_3} - \frac{q_1}{t_1},\\
	&b_3 = t_2p_2 - (q_3p_3-\rho_3)\frac{t_2p_2+q_3p_1}{\rho_2-\rho_3},\quad
	c_3 = -\frac{q_2}{t_2},\\
	&b_4 = (q_3p_3-\rho_3)\frac{t_2p_2-q_3p_1}{\rho_2-\rho_3},\quad
	c_4 = - \frac{q_2}{t_2} + \frac{\rho_2-\rho_3}{t_2p_2-q_3p_1}.
\end{split}\]
Although the components of the matrices $\widetilde{A}_1,\ldots,\widetilde{A}_4$ are rational in $(q_j,p_j)$, the Hamiltonians $\widetilde{H}_1,\widetilde{H}_2$ turn out to be polynomials in $(q_j,p_j)$.
\begin{thm}
The dependent variables $q_j,p_j$ $(j=1,2,3)$ satisfy the system $\mathcal{H}^{21,21,21,21,111}$.
\end{thm}

\subsection{Spectral type $31,31,22,22,22$}

In this case, we consider a Fuchsian system
\begin{equation}\label{Eq:Fuchs_31_31_22_22_22}
	\frac{\partial}{\partial x}Y(x) = \left(\frac{A_1}{x-t_1}+\frac{A_2}{x-t_2}+\frac{A_3}{x-1}+\frac{A_4}{x}\right)Y(x),
\end{equation}
with a Riemann scheme
\[
	\left\{\begin{array}{ccccc}
		x=t_1 & x=t_2 & x=1 & x=0 & x=\infty \\
		\theta_1 & \theta_2 & \theta_3 & \theta_4 & \rho_1 \\
		0 & 0 & \theta_3 & \theta_4 & \rho_1 \\
		0 & 0 & 0 & 0 & \rho_2 \\
		0 & 0 & 0 & 0 & \rho_2
	\end{array}\right\}.
\]
Note that a Fuchsian relation $\theta_1+\theta_2+2\theta_3+2\theta_4+2\rho_1+2\rho_2=0$ is satisfied.

In a similar manner as Section \ref{Sec:Derivation_21_21_21_21_111}, the residue matrices are transformed to
\[\begin{split}
	&\widetilde{A}_1 = \begin{pmatrix}1\\0\\b_1\\b_2\end{pmatrix}\begin{pmatrix}\theta_1-b_1c_1-b_2c_2&a_1&c_1&c_2\end{pmatrix},\quad
	\widetilde{A}_2 = \begin{pmatrix}0\\1\\a_2\\b_3\end{pmatrix}\begin{pmatrix}a_3&\theta_2-a_2-b_3c_3&1&c_3\end{pmatrix},\\
	&\widetilde{A}_3 = \begin{pmatrix}1&0\\0&1\\a_4&a_5\\a_6&a_7\end{pmatrix}\begin{pmatrix}\theta_3-a_4&-a_5&1&0\\-a_6&\theta_3-a_7&0&1\end{pmatrix},\quad
	\widetilde{A}_4 = \begin{pmatrix}1&0\\0&1\\0&0\\0&0\end{pmatrix}\begin{pmatrix}\theta_4&0&a_8&a_9\\0&\theta_4&a_{10}&a_{11}\end{pmatrix},
\end{split}\]
and
\begin{equation}\label{Eq:Residue_31_31_22_22_22}
	\widetilde{A}_{\infty} := -\sum_{i=1}^{4}\widetilde{A}_i = \begin{pmatrix}\rho_1&0&0&0\\0&\rho_1&0&0\\a_{12}&a_{13}&\rho_2&0\\a_{14}&a_{15}&0&\rho_2\end{pmatrix}.
\end{equation}
By using \eqref{Eq:Residue_31_31_22_22_22}, we can show that the variables $a_1,\ldots.a_{15}$ are polynomials in $(b_j,c_j)$; we do not give their explicit formulas here.
Then the dependent variables $b_j,c_j$ $(j=1,2,3)$ satisfy a Hamiltonian system
\begin{equation}\begin{split}\label{Eq:Sch_Ham_31_31_22_22_22}
	&\frac{\partial b_j}{\partial t_i} = \{H_i,b_j\},\quad
	\frac{\partial c_j}{\partial t_i} = \{H_i,c_j\}\quad (i=1,2),\\
	&H_1 = \frac{\mathrm{tr}\widetilde{A}_1\widetilde{A}_2}{t_1-t_2} + \frac{\mathrm{tr}\widetilde{A}_1\widetilde{A}_3}{t_1-1} + \frac{\mathrm{tr}\widetilde{A}_1\widetilde{A}_4}{t_1},\quad
	H_2 = \frac{\mathrm{tr}\widetilde{A}_2\widetilde{A}_1}{t_2-t_1} + \frac{\mathrm{tr}\widetilde{A}_2\widetilde{A}_3}{t_2-1} + \frac{\mathrm{tr}\widetilde{A}_2\widetilde{A}_4}{t_2},
\end{split}\end{equation}
with a symplectic form $\omega=\sum_{j=1}^{3}db_j\wedge dc_j$.

Under the system \eqref{Eq:Sch_Ham_31_31_22_22_22}, we consider a canonical transformation
\[\begin{split}
	&\frac{q_1}{t_1} = -c_1,\quad
	t_1p_1 = b_1+b_2c_3,\quad
	\frac{q_2}{t_2} = -c_3,\quad
	t_2p_2 = b_3+b_2c_1,\\
	&\frac{q_3}{t_1t_2} = -c_2 + c_1c_3,\quad
	t_1t_2p_3 = b_2.
\end{split}\]
Then, by a direct computation, we arrive at
\begin{thm}
The dependent variables $q_j,p_j$ $(j=1,2,3)$ satisfy the system $\mathcal{H}^{31,31,22,22,22}$.
\end{thm}

\subsection{Spectral type $21,111,111,111$}

In this case, we consider a Fuchsian system
\[
	\frac{\partial}{\partial x}Y(x) = \left(\frac{A_1}{x-t}+\frac{A_2}{x-1}+\frac{A_3}{x}\right)Y(x),
\]
with a Riemann scheme
\[
	\left\{\begin{array}{cccc}
		x=t & x=1 & x=0 & x=\infty \\
		\theta_1 & \theta_{2,1} & \theta_{3,1} & \rho_1 \\
		0 & \theta_{2,2} & \theta_{3,2} & \rho_2 \\
		0 & 0 & 0 & \rho_3
	\end{array}\right\}.
\]
Note that a Fuchsian relation $\theta_1+\theta_{2,1}+\theta_{2,2}+\theta_{3,1}+\theta_{3,2}+\rho_1+\rho_2+\rho_3=0$ is satisfied.

In a similar manner as Section \ref{Sec:Derivation_21_21_21_21_111}, the residue matrices are transformed to
\[\begin{split}
	&\widetilde{A}_1= \begin{pmatrix}1\\b_1\\b_2\end{pmatrix}\begin{pmatrix}\theta_1-b_1c_1-b_2c_2&c_1&c_2\end{pmatrix},\quad
	\widetilde{A}_2 = \begin{pmatrix}1&1\\b_3&a_1\\b_4&a_2\end{pmatrix}\begin{pmatrix}a_3&c_3&c_4\\a_4&1&1\end{pmatrix},\\
	&\widetilde{A}_3 = \begin{pmatrix}1&0\\0&1\\0&0\end{pmatrix}\begin{pmatrix}\theta_{3,1}&a_5&a_6\\0&\theta_{3,2}&a_7\end{pmatrix},
\end{split}\]
where
\begin{equation}\label{Eq:A2_21_111_111_111}
	\begin{pmatrix}a_3&c_3&c_4\\a_4&1&1\end{pmatrix}\begin{pmatrix}1&1\\b_3&a_1\\b_4&a_2\end{pmatrix} = \begin{pmatrix}\theta_{2,1}&0\\0&\theta_{2,2}\end{pmatrix},
\end{equation}
and
\begin{equation}\label{Eq:Residue_21_111_111_111}
	\widetilde{A}_{\infty} := -\sum_{i=1}^{3}\widetilde{A}_i = \begin{pmatrix}\rho_1&0&0\\a_8&\rho_2&0\\a_9&a_{10}&\rho_3\end{pmatrix}.
\end{equation}
By using \eqref{Eq:A2_21_111_111_111} and \eqref{Eq:Residue_21_111_111_111}, we can show that the variables $a_1,\ldots.a_{10}$ are rational in $(b_j,c_j)$; we do not give their explicit formulas here.
Furthermore, we obtain
\begin{equation}\begin{split}\label{Eq:Residue_21_111_111_111_Comp}
	&b_1c_1 + b_2c_2 + b_3(c_3+1)  + b_4(c_4+1) + \theta_{2,2} + \theta_{3,2} + \rho_2 + \rho_3 = 0,\\
	&\{b_2c_2+b_4(c_4+1)+\rho_3\}(c_3-c_4) + \theta_{2,2}c_3 + \theta_{2,1}= 0.
\end{split}\end{equation}
Then the dependent variables $b_j,c_j$ $(j=1,\ldots,4)$ satisfy a Hamiltonian system\begin{equation}\label{Eq:Sch_Ham_21_111_111_111}
	\frac{\partial b_j}{\partial t} = \{H,b_j\},\quad
	\frac{\partial c_j}{\partial t} = \{H,c_j\},\quad
	H = \frac{\mathrm{tr}\widetilde{A}_1\widetilde{A}_2}{t-1} + \frac{\mathrm{tr}\widetilde{A}_1\widetilde{A}_3}{t},
\end{equation}
with a symplectic form
\begin{equation}\label{Eq:Simplectic_21_111_111_111}
	\omega = \sum_{j=1}^{4}db_j\wedge dc_j,
\end{equation}
and the relation \eqref{Eq:Residue_21_111_111_111_Comp}.
Note that the Hamiltonian $H$ turns out to be a polynomial in $(b_j,c_j)$, although the components of the matrices $\widetilde{A}_1,\widetilde{A}_2,\widetilde{A}_3$ are rational.

We reduce the Hamiltonian system \eqref{Eq:Sch_Ham_21_111_111_111} to the one of sixth order.
Substituting the first relation of \eqref{Eq:Residue_21_111_111_111_Comp} to \eqref{Eq:Simplectic_21_111_111_111}, we obtain
\[\begin{split}
	\omega &= db_1\wedge dc_1 + db_2\wedge dc_2 + db_3\wedge d(c_3+1) + db_4\wedge d(c_4+1)\\
	&= db_1\wedge dc_1 + db_2\wedge dc_2 + db_3\wedge d(c_3+1) - d\frac{b_1c_1+b_2c_2+b_3(c_3+1)}{c_4+1}\wedge d(c_4+1)\\
	&= db_1(c_4+1)\wedge d\frac{c_1}{c_4+1} + db_2(c_4+1)\wedge d\frac{c_2}{c_4+1} + db_3(c_4+1)\wedge d\frac{c_3+1}{c_4+1}.
\end{split}\]
Hence we can take a six-dimensional canonical coordinate system by
\[\begin{split}
	&\frac{q_1}{t} = -\frac{c_1}{c_4+1},\quad
	tp_1 = b_1(c_4+1),\quad
	\frac{q_2}{t} = -\frac{c_2}{c_4+1},\quad
	tp_2 = b_2(c_4+1),\\
	&q_3-1 = -\frac{c_3+1}{c_4+1},\quad
	p_3 = b_3(c_4+1).
\end{split}\]
Note that the variables $b_j,c_j$ $(j=1,\ldots,4)$ are rational in $(q_j,p_j)$; we do not give their explicit formulas here.
Then, in a similar manner as Section \ref{Sec:Derivation_21_21_21_21_111}, we arrive at
\begin{thm}
The dependent variables $q_j,p_j$ $(j=1,2,3)$ satisfy the system $\mathcal{H}^{21,111,111,111}$.
\end{thm}

\subsection{Spectral type $31,22,211,1111$}

In this case, we consider a Fuchsian system
\begin{equation}\label{Eq:Fuchs_31_22_211_1111}
	\frac{\partial}{\partial x}Y(x) = \left(\frac{A_1}{x-t}+\frac{A_2}{x-1}+\frac{A_3}{x}\right)Y(x),
\end{equation}
with a Riemann scheme
\[
	\left\{\begin{array}{cccc}
		x=t & x=1 & x=0 & x=\infty \\
		\theta_1 & \theta_2 & \theta_{3,1} & \rho_1 \\
		0 & \theta_2 & \theta_{3,2} & \rho_2 \\
		0 & 0 & 0 & \rho_3 \\
		0 & 0 & 0 & \rho_4 \\
	\end{array}\right\}.
\]
Note that a Fuchsian relation $\theta_1+2\theta_2+\theta_{3,1}+\theta_{3,2}+\rho_1+\rho_2+\rho_3+\rho_4=0$ is satisfied.

In a similar manner as Section \ref{Sec:Derivation_21_21_21_21_111}, the residue matrices are transformed to
\[\begin{split}
	&\widetilde{A}_1= \begin{pmatrix}1\\b_1\\b_2\\b_3\end{pmatrix}\begin{pmatrix}\theta_1-b_1c_1-b_2c_2-b_3c_3&c_1&c_2&c_3\end{pmatrix},\\
	&\widetilde{A}_2 = \begin{pmatrix}1&0\\0&1\\a_1&a_2\\a_3&b_4\end{pmatrix}\begin{pmatrix}\theta_2-a_1-a_3&-a_2-b_4&1&1\\-a_1-c_4a_3&\theta_2-a_2-b_4c_4&1&c_4\end{pmatrix},\quad
	\widetilde{A}_3 = \begin{pmatrix}1&0\\0&1\\0&0\\0&0\end{pmatrix}\begin{pmatrix}\theta_{3,1}&a_4&a_5&a_6\\0&\theta_{3,2}&a_7&a_8\end{pmatrix},
\end{split}\]
and
\begin{equation}\label{Eq:Residue_31_22_211_1111}
	\widetilde{A}_{\infty} := -\sum_{i=1}^{3}\widetilde{A}_i = \begin{pmatrix}\rho_1&0&0&0\\a_9&\rho_2&0&0\\a_{10}&a_{11}&\rho_3&0\\a_{12}&a_{13}&0&\rho_4\end{pmatrix}.
\end{equation}
By using \eqref{Eq:Residue_31_22_211_1111}, we can show that the variables $a_1,\ldots.a_{13}$ are polynomials in $(b_j,c_j)$; we do not give their explicit formulas here.
Furthermore, we obtain
\begin{equation}\begin{split}\label{Eq:Residue_31_22_211_1111_Comp}
	&(b_1c_1-b_4c_4+\theta_2+\theta_{3,2}+\rho_2)(c_4-1) - b_2(c_2-c_3) - \rho_3 = 0,\\
	&b_3(c_3-c_2) + b_4(c_4-1) + \rho_4 = 0.
\end{split}\end{equation}
Then the dependent variables $b_j,c_j$ $(j=1,\ldots,4)$ satisfy a Hamiltonian system
\begin{equation}\label{Eq:Sch_Ham_31_22_211_1111}
	\frac{\partial b_j}{\partial t} = \{H,b_j\},\quad
	\frac{\partial c_j}{\partial t} = \{H,c_j\},\quad
	H = \frac{\mathrm{tr}\widetilde{A}_1\widetilde{A}_2}{t-1} + \frac{\mathrm{tr}\widetilde{A}_1\widetilde{A}_3}{t},
\end{equation}
with a symplectic form
\begin{equation}\label{Eq:Simplectic_31_22_211_1111}
	\omega = \sum_{j=1}^{4}db_j\wedge dc_j,
\end{equation}
and the relation \eqref{Eq:Residue_31_22_211_1111_Comp}.

We reduce the Hamiltonian system \eqref{Eq:Sch_Ham_31_22_211_1111} to the one of sixth order.
Substituting the second relation of \eqref{Eq:Residue_31_22_211_1111_Comp} to \eqref{Eq:Simplectic_31_22_211_1111}, we obtain
\[\begin{split}
	\omega &= db_1\wedge dc_1 + d(b_2+b_3)\wedge dc_2 + db_3\wedge d(c_3-c_2) + db_4\wedge d(c_4-1)\\
	&= db_1\wedge dc_1 + d(b_2+b_3)\wedge dc_2 + db_3\wedge d(c_3-c_2) - d\frac{b_3(c_3-c_2)}{c_4-1}\wedge d(c_4-1)\\
	&= db_1\wedge dc_1 + d(b_2+b_3)\wedge dc_2 + db_3(c_4-1)\wedge d\frac{c_3-c_2}{c_4-1}.
\end{split}\]
Hence we can take a six-dimensional canonical coordinate system by
\[
	\lambda_1 = b_1,\quad
	\mu_1 = c_1,\quad
	\lambda_2 = -c_2,\quad
	\mu_2 = b_2+b_3,\quad
	\lambda_3 = -\frac{c_3-c_2}{c_4-1},\quad
	\mu_3 = b_3(c_4-1).
\]
Furthermore, we consider a canonical transformation
\[\begin{split}
	&\frac{q_1}{t} = \frac{1}{\lambda_1},\quad
	tp_1 = -\lambda_1(\lambda_1\mu_1-\lambda_3\mu_3+\theta_2+\theta_{3,2}+\rho_2+\rho_4),\\
	&\frac{q_2}{t} = \lambda_2,\quad
	tp_2 = \mu_2,\quad
	\frac{q_3}{t} = \lambda_1\lambda_3,\quad
	tp_3 = \frac{\mu_3}{\lambda_1}.
\end{split}\]
Note that the variables $b_j,c_j$ $(j=1,\ldots,4)$ are rational in $(q_j,p_j)$; we do not give their explicit formulas here.
Then, in a similar manner as Section \ref{Sec:Derivation_21_21_21_21_111}, we arrive at
\begin{thm}
The dependent variables $q_j,p_j$ $(j=1,2,3)$ satisfy the system $\mathcal{H}^{31,22,211,1111}$.
\end{thm}

\subsection{Spectral type $22,22,211,211$}\label{Sec:Derivation_22_22_211_211}

In this case, we consider a Fuchsian system
\[
	\frac{\partial}{\partial x}Y(x) = \left(\frac{A_1}{x-t}+\frac{A_2}{x-1}+\frac{A_3}{x}\right)Y(x),
\]
with a Riemann scheme
\[
	\left\{\begin{array}{cccc}
		x=t & x=1 & x=0 & x=\infty \\
		\theta_1 & \theta_2 & \theta_{3,1} & \rho_1 \\
		\theta_1 & \theta_2 & \theta_{3,2} & \rho_2 \\
		0 & 0 & 0 & \rho_3 \\
		0 & 0 & 0 & \rho_3 \\
	\end{array}\right\}.
\]
Note that a Fuchsian relation $2\theta_1+2\theta_2+\theta_{3,1}+\theta_{3,2}+\rho_1+\rho_2+2\rho_3=0$ is satisfied.

In a similar manner as Section \ref{Sec:Derivation_21_21_21_21_111}, the residue matrices are transformed to
\[\begin{split}
	&\widetilde{A}_1 = \begin{pmatrix}I_2\\\widetilde{B}_1\end{pmatrix}\begin{pmatrix}\theta_1I_2-\widetilde{C}_1\widetilde{B}_1&\widetilde{C}_1\end{pmatrix},\quad
	\widetilde{B}_1 = \begin{pmatrix}b_1&b_2\\a_1&b_3\end{pmatrix},\quad
	\widetilde{C}_1 = \begin{pmatrix}c_1&1\\c_2&c_3\end{pmatrix},\\
	&\widetilde{A}_2 = \begin{pmatrix}I_2\\\widetilde{B}_2\end{pmatrix}\begin{pmatrix}\theta_2I_2-\widetilde{B}_2&I_2\end{pmatrix},\quad
	\widetilde{B}_2 = \begin{pmatrix}a_2&a_3\\a_4&a_5\end{pmatrix},\\
	&\widetilde{A}_3 = \begin{pmatrix}I_2\\O\end{pmatrix}\begin{pmatrix}\widetilde{C}_{3,1}&\widetilde{C}_{3,2}\end{pmatrix},\quad
	\widetilde{C}_{3,1} = \begin{pmatrix}\theta_{3,1}&a_6\\0&\theta_{3,2}\end{pmatrix},\quad
	\widetilde{C}_{3,2} = \begin{pmatrix}a_7&a_8\\a_9&a_{10}\end{pmatrix}.
\end{split}\]
and
\begin{equation}\label{Eq:Residue_22_22_211_211}
	\widetilde{A}_{\infty} := -\sum_{i=1}^{3}\widetilde{A}_i = \begin{pmatrix}\rho_1&0&0&0\\a_{11}&\rho_2&0&0\\a_{12}&a_{13}&\rho_3&0\\a_{14}&a_{15}&0&\rho_3\end{pmatrix}.
\end{equation}
By using \eqref{Eq:Residue_22_22_211_211}, we can show that the variables $a_1,\ldots.a_{15}$ are polynomials in $(b_j,c_j)$; we do not give their explicit formulas here.
Then the dependent variables $b_j,c_j$ $(j=1,2,3)$ satisfy a Hamiltonian system
\begin{equation}\label{Eq:Sch_Ham_22_22_211_211}
	\frac{\partial b_j}{\partial t} = \{H,b_j\},\quad
	\frac{\partial c_j}{\partial t} = \{H,c_j\},\quad
	H = \frac{\mathrm{tr}\widetilde{A}_1\widetilde{A}_2}{t-1} + \frac{\mathrm{tr}\widetilde{A}_1\widetilde{A}_3}{t},
\end{equation}
with a symplectic form $\omega=\sum_{j=1}^{3}db_j\wedge dc_j$.

The system \eqref{Eq:Sch_Ham_22_22_211_211} is transformed into the one $\mathcal{H}^{22,22,211,211}$ as follows.
The Hamiltonian $H$ is described as
\[\begin{split}
	H &= \frac{1}{t-1}\mathrm{tr}(\theta_2I_2-\widetilde{B}_2+\widetilde{B}_1)(\theta_1I_2-\widetilde{C}_1\widetilde{B}_1+\widetilde{C}_1\widetilde{B}_2) + \frac{1}{t}\mathrm{tr}(\widetilde{C}_{3,1}+\widetilde{C}_{3,2}\widetilde{B}_1)(\theta_1I_2-\widetilde{C}_1\widetilde{B}_1)\\
	&= \frac{1}{t-1}\mathrm{tr}\{(\theta_2+\rho_3)I_2+\widetilde{B}_1(I_2+\widetilde{C}_1)\}\{\theta_1I_2-\rho_3\widetilde{C}_1-\widetilde{C}_1\widetilde{B}_1(I_2+\widetilde{C}_1)\}\\
	&\quad + \frac{1}{t}\mathrm{tr}\{\widetilde{C}_{3,1}-(I_2+\widetilde{C}_1)\widetilde{B}_1\}(\theta_1I_2-\widetilde{C}_1\widetilde{B}_1)\\
	&= \frac{1}{t-1}\mathrm{tr}\theta_1(\theta_2+\rho_3)I_2 - \frac{1}{t-1}\mathrm{tr}\rho_3(\theta_2+\rho_3)\widetilde{C}_1 - \frac{1}{t-1}\mathrm{tr}(\theta_2+2\rho_3)\widetilde{B}_1\widetilde{C}_1(I_2+\widetilde{C}_1)\\
	&\quad + \frac{1}{t(t-1)}\mathrm{tr}\theta_1\widetilde{B}_1(I_2+\widetilde{C}_1) - \frac{1}{t(t-1)}\mathrm{tr}\widetilde{B}_1(I_2+\widetilde{C}_1)\widetilde{B}_1\widetilde{C}_1(I_2+t\widetilde{C}_1) - \frac{1}{t}\mathrm{tr}\widetilde{C}_{3,1}(\widetilde{C}_1\widetilde{B}_1-\theta_1I_2).
\end{split}\]
Here we set
\[\begin{split}
	q_1 = c_1,\quad
	p_1 = -b_1,\quad
	q_2 = c_2,\quad
	p_2 = -b_2,\quad
	q_3 = c_3,\quad
	p_3 = -b_3,
\end{split}\]
and
\[
	P = \frac{1}{t}\widetilde{B}_1,\quad
	Q = -t\widetilde{C}_1,\quad
	R = -\theta_{3,1}I_2+\widetilde{C}_{3,1}.
\]
Then it is easy to verify that the Hamiltonian $H$ is just equivalent to the one $H^{22,22,211,211}$, which was given in Section \ref{Sec:Hamiltonian}.
Note that
\[
	a_1 = q_2p_2 - \theta_1 - \theta_2 - \theta_{3,2} - \rho_2 - \rho_3,\quad
	a_6 = -(q_1-q_3)p_2 + p_1 - p_3.
\]
\begin{thm}
The dependent variables $q_j,p_j$ $(j=1,2,3)$ satisfy the system $\mathcal{H}^{22,22,211,211}$.
\end{thm}

\subsection{Spectral type $22,22,22,1111$}

In this case, we consider a Fuchsian system
\[
	\frac{\partial}{\partial x}Y(x) = \left(\frac{A_1}{x-t}+\frac{A_2}{x-1}+\frac{A_3}{x}\right)Y(x),
\]
with a Riemann scheme
\[
	\left\{\begin{array}{cccc}
		x=t & x=1 & x=0 & x=\infty \\
		\theta_1 & \theta_2 & \theta_3 & \rho_1 \\
		\theta_1 & \theta_2 & \theta_3 & \rho_2 \\
		0 & 0 & 0 & \rho_3 \\
		0 & 0 & 0 & \rho_4 \\
	\end{array}\right\}.
\]
Note that a Fuchsian relation $2\theta_1+2\theta_2+2\theta_3+\rho_1+\rho_2+\rho_3+\rho_4=0$ is satisfied.

In a similar manner as Section \ref{Sec:Derivation_21_21_21_21_111}, the residue matrices are transformed to
\[\begin{split}
	&\widetilde{A}_1 = \begin{pmatrix}I_2\\\widetilde{B}_1\end{pmatrix}\begin{pmatrix}\theta_1I_2-\widetilde{C}_1\widetilde{B}_1&\widetilde{C}_1\end{pmatrix},\quad
	\widetilde{B}_1 = \begin{pmatrix}b_1&b_2\\b_3&b_4\end{pmatrix},\quad
	\widetilde{C}_1 = \begin{pmatrix}c_1&c_3\\c_2&c_4\end{pmatrix},\\
	&\widetilde{A}_2 = \begin{pmatrix}I_2\\\widetilde{B}_2\end{pmatrix}\begin{pmatrix}\theta_2I_2-\widetilde{C}_2\widetilde{B}_2&\widetilde{C}_2\end{pmatrix},\quad
	\widetilde{B}_2 = \begin{pmatrix}a_1&b_5\\a_2&a_3\end{pmatrix},\quad
	\widetilde{C}_2 = \begin{pmatrix}1&1\\c_5&1\end{pmatrix},\\
	&\widetilde{A}_3 = \begin{pmatrix}I_2\\O\end{pmatrix}\begin{pmatrix}\theta_3I_2&\widetilde{C}_3\end{pmatrix},\quad
	\widetilde{C}_3 = \begin{pmatrix}a_4&a_5\\a_6&a_7\end{pmatrix},
\end{split}\]
and
\begin{equation}\label{Eq:Residue_22_22_22_1111}
	\widetilde{A}_{\infty} := -\sum_{i=1}^{3}\widetilde{A}_i = \begin{pmatrix}\rho_1&0&0&0\\0&\rho_2&0&0\\a_8&a_9&\rho_3&0\\a_{10}&a_{11}&0&\rho_4\end{pmatrix}.
\end{equation}
By using \eqref{Eq:Residue_22_22_22_1111}, we can show that the variables $a_1,\ldots.a_{11}$ are polynomials in $(b_j,c_j)$; we do not give their explicit formulas here.
Furthermore, we obtain
\begin{equation}\begin{split}\label{Eq:Residue_22_22_22_1111_Comp}
	&b_2(c_2-c_1) + b_4(c_4-c_3) + b_5(c_5-1) - \theta_1 - \theta_2 - \theta_3 - \rho_2 = 0,\\
	&b_1(c_1c_5-c_2) + b_2c_2(c_5-1) + b_3(c_3-c_4) + b_5c_5(c_5-1) - \rho_3c_5 + \theta_1 + \theta_2 + \theta_3 + \rho_2 - \rho_4 = 0,\\
	&b_1(c_1-c_3) + b_2(c_2-c_4) + b_5(c_5-1) - \rho_3 = 0,\\
	&b_2c_2(c_5-1) + b_3(c_3-c_1) + b_4(c_4c_5-c_2) + b_5c_5(c_5-1) - (\theta_1+\theta_2+\theta_3+\rho_2)(c_5-1) - \rho_4 = 0.
\end{split}\end{equation}
Then the dependent variables $b_j,c_j$ $(j=1,\ldots,5)$ satisfy a Hamiltonian system
\begin{equation}\label{Eq:Sch_Ham_22_22_22_1111}
	\frac{\partial b_j}{\partial t} = \{H,b_j\},\quad
	\frac{\partial c_j}{\partial t} = \{H,c_j\},\quad
	H = \frac{\mathrm{tr}\widetilde{A}_1\widetilde{A}_2}{t-1} + \frac{\mathrm{tr}\widetilde{A}_1\widetilde{A}_3}{t},
\end{equation}
with a symplectic form
\begin{equation}\label{Eq:Simplectic_22_22_22_1111}
	\omega = \sum_{j=1}^{5}db_j\wedge dc_j,
\end{equation}
and the relation \eqref{Eq:Residue_22_22_22_1111_Comp}.

We derive a six-dimensional canonical coordinate system in advance.
The first and third relation of \eqref{Eq:Residue_22_22_22_1111_Comp} are rewritten to
\[\begin{split}
	&(b_1+b_2)(c_1-c_3) - (b_4+b_2)(c_4-c_3) + \theta_1 + \theta_2 + \theta_3 + \rho_2 - \rho_3 = 0,\\
	&(b_1+b_2)(c_1-c_3) + b_2(c_2-c_4-c_1+c_3) + b_5(c_5-1) - \rho_3 = 0.
\end{split}\]
Substituting them to \eqref{Eq:Simplectic_22_22_22_1111}, we obtain
\[\begin{split}
	\omega &= d(b_1+b_2)\wedge d(c_1-c_3) + db_2\wedge d(c_2-c_4-c_1+c_3) + d(b_3+b_4+b_1+b_2)\wedge dc_3\\
	&\quad + d(b_4+b_2)\wedge d(c_4-c_3) + db_5\wedge d(c_5-1)\\
	&= d(b_1+b_2)\wedge d(c_1-c_3) + db_2\wedge d(c_2-c_4-c_1+c_3) + d(b_3+b_4+b_1+b_2)\wedge dc_3\\
	&\quad + d\frac{(b_1+b_2)(c_1-c_3)}{c_4-c_3}\wedge d(c_4-c_3)\\
	&\quad - d\{(b_1+b_2)(c_1-c_3)+b_2(c_2-c_4-c_1+c_3)\}\wedge d(c_5-1)\\
	&= d\frac{(b_1+b_2)(c_5-1)}{c_4-c_3}\wedge d\frac{(c_1-c_3)(c_4-c_3)}{c_5-1} + db_2(c_5-1)\wedge d\frac{c_2-c_4-c_1+c_3}{c_5-1}\\
	&\quad + d(b_3+b_4+b_1+b_2)\wedge dc_3.
\end{split}\]
Hence we can take
\begin{equation}\begin{split}\label{Eq:Cano_Coor_22_22_22_1111}
	&q_1 = c_3,\quad
	p_1 = -b_3 - b_4 - b_1 - b_2,\\
	&q_2 = \frac{(c_1-c_3)(c_4-c_3)}{c_5-1},\quad
	p_2 + \frac{\theta_1+\theta_2+\theta_3+\rho_2+\rho_3}{q_2} = -\frac{(b_1+b_2)(c_5-1)}{c_4-c_3},\\
	&q_3 = \frac{c_2-c_4-c_1+c_3}{c_5-1},\quad
	p_3 = -b_2(c_5-1).
\end{split}\end{equation}

The system \eqref{Eq:Sch_Ham_22_22_22_1111} is transformed into the one $\mathcal{H}^{22,22,22,1111}$ as follows.
We set
\[
	P = \frac{1}{t}\Gamma^{-1}\widetilde{B}_1\widetilde{C}_2\Gamma,\quad
	Q = -t\Gamma^{-1}\widetilde{C}_2^{-1}\widetilde{B}_1\Gamma,\quad
	R = \Gamma^{-1}\begin{pmatrix}0&0\\0&-\rho_3+\rho_4\end{pmatrix}\Gamma,
\]
where
\[
	\Gamma = \begin{pmatrix}0&-1\\c_1-c_3&1\end{pmatrix}.
\]
Then, by using \eqref{Eq:Residue_22_22_22_1111_Comp} and \eqref{Eq:Cano_Coor_22_22_22_1111}, we can show that the components of the matrices $P,Q,R$ are polynomials in $(q_j,p_j)$; we do not give their explicit formulas here.
And, in a similar manner as Section \ref{Sec:Derivation_22_22_211_211}, we arrive at
\begin{thm}
The dependent variables $q_j,p_j$ $(j=1,2,3)$ satisfy the system $\mathcal{H}^{22,22,22,1111}$.
\end{thm}

\subsection{Spectral type $42,33,33,222$}

In this case, we consider a Fuchsian system
\[
	\frac{\partial}{\partial x}Y(x) = \left(\frac{A_1}{x-t}+\frac{A_2}{x-1}+\frac{A_3}{x}\right)Y(x),
\]
with a Riemann scheme
\[
	\left\{\begin{array}{cccc}
		x=t & x=1 & x=0 & x=\infty \\
		\theta_1 & \theta_2 & \theta_3 & \rho_1 \\
		\theta_1 & \theta_2 & \theta_3 & \rho_1 \\
		\theta_1 & \theta_2 & 0 & \rho_2 \\
		0 & 0 & 0 & \rho_2 \\
		0 & 0 & 0 & \rho_3 \\
		0 & 0 & 0 & \rho_3 \\
	\end{array}\right\}.
\]
Note that a Fuchsian relation $3\theta_1+3\theta_2+2\theta_3+2\rho_1+2\rho_2+2\rho_3=0$ is satisfied.

In a similar manner as Section \ref{Sec:Derivation_21_21_21_21_111}, the residue matrices are transformed to
\[\begin{split}
	&\widetilde{A}_1 = \begin{pmatrix}I_3\\\widetilde{B}_1\end{pmatrix}\begin{pmatrix}\theta_1I_3-\widetilde{C}_1\widetilde{B}_1&\widetilde{C}_1\end{pmatrix},\quad
	\widetilde{B}_1 = \begin{pmatrix}0&b_1&0\\b_2&a_1&b_3\\b_4&b_2&b_5\end{pmatrix},\quad
	\widetilde{C}_1 = \begin{pmatrix}a_3&c_2&c_4\\c_1&0&1\\a_4&c_3&c_5\end{pmatrix},\\
	&\widetilde{A}_2 = \begin{pmatrix}I_3\\\widetilde{B}_2\end{pmatrix}\begin{pmatrix}\theta_2I_3-\widetilde{C}_2\widetilde{B}_2&\widetilde{C}_2\end{pmatrix},\quad
	\widetilde{B}_2 = \begin{pmatrix}a_5&a_6&1\\a_7&a_8&b_6\\a_9&a_{10}&b_7\end{pmatrix},\quad
	\widetilde{C}_2 = \begin{pmatrix}1&0&1\\0&1&0\\a_{11}&c_6&c_7\end{pmatrix},\\
	&\widetilde{A}_3 = \begin{pmatrix}1&0\\0&1\\0&0\\0&0\\0&0\\0&0\end{pmatrix}\begin{pmatrix}\theta_3&0&a_{12}&a_{13}&a_{14}&a_{15}\\0&\theta_3&a_{16}&a_{17}&a_{18}&a_{19}\end{pmatrix},
\end{split}\]
and
\begin{equation}\label{Eq:Residue_42_33_33_222}
	\widetilde{A}_{\infty} := -\sum_{i=1}^{3}\widetilde{A}_i = \begin{pmatrix}\rho_1&0&0&0&0&0\\0&\rho_1&0&0&0&0\\a_{20}&a_{21}&\rho_2&0&0&0\\a_{22}&a_{23}&0&\rho_2&0&0\\a_{24}&a_{25}&0&0&\rho_3&0\\a_{26}&a_{27}&0&0&0&\rho_3\end{pmatrix}.
\end{equation}
Now we can find 35 relations in \eqref{Eq:Residue_42_33_33_222}.
Among them, the relations derived from the following matrix components are used to determine the variables $a_{12},\ldots,a_{27}$ as polynomials in the other variables:
\[\begin{array}{llllllll}
	(1,3), & (1,4), & (1,5), & (1,6), & (2,3), & (2,4), & (2,5), & (2,6), \\
	(3,1), & (3,2), & (4,1), & (4,2), & (5,1), & (5,2), & (6,1), & (6,2).
\end{array}\]
And the following ones are used to determine the variables $a_1,\ldots,a_{11}$ as rational expressions in $(b_j,c_j)$:
\[\begin{array}{lllllllllll}
	(1,2), & (2,1), & (2,2), & (3,3), & (3,4), & (4,4), & (4,5), & (5,5), & (5,6), & (6,5), & (6,6).
\end{array}\]
We do not give their explicit formulas here.
Then the dependent variables $b_j,c_j$ $(j=1,\ldots,7)$ satisfy a rational Hamiltonian system
\begin{equation}\label{Eq:Sch_Ham_42_33_33_222}
	\frac{\partial b_j}{\partial t} = \{H,b_j\},\quad
	\frac{\partial c_j}{\partial t} = \{H,c_j\},\quad
	H = \frac{\mathrm{tr}\widetilde{A}_1\widetilde{A}_2}{t-1} + \frac{\mathrm{tr}\widetilde{A}_1\widetilde{A}_3}{t},
\end{equation}
with a symplectic form
\begin{equation}\label{Eq:Simplectic_42_33_33_222_7th}
	\omega = \sum_{j=1}^{7}db_j\wedge dc_j.
\end{equation}
Furthermore, we have 8 relations which is derived from the matrix components
\[\begin{array}{llllllll}
	(3,5), & (3,6), & (4,3), & (4,6), & (5,3), & (5,4), & (6,3), & (6,4).
\end{array}\]
In order to derive the Hamiltonian system of sixth order, we use the first four relations, whose explicit formulas are given as
\begin{equation}\begin{split}\label{Eq:Residue_42_33_33_222_7th}
	&c_3 + c_6 = 0,\\
	&c_5 + c_7 = 0,\\
	&b_1c_1 - b_1c_5 + \theta_2 + \rho_2 = 0,\\
	&b_1(c_1-1) - b_3c_3 - b_5c_5 - b_6c_6 - b_7(c_7+1) + \theta_1 + \theta_2 + 2\rho_2 = 0.
\end{split}\end{equation}

We reduce the Hamiltonian system \eqref{Eq:Sch_Ham_42_33_33_222} to the one of sixth order.
Substituting the first and second relation of \eqref{Eq:Residue_42_33_33_222_7th} to \eqref{Eq:Simplectic_42_33_33_222_7th}, we obtain
\[\begin{split}
	\omega &= db_1\wedge d(c_1-1) + db_2\wedge dc_2 + d(b_3-b_6)\wedge dc_3 + db_4\wedge dc_4 + d(b_5-b_7-1)\wedge dc_5.
\end{split}\]
Hence we can take
\[\begin{split}
	&\lambda_1=-c_1+1,\quad \mu_1=b_1,\quad \lambda_2=-c_2,\quad \mu_2=b_2,\quad \lambda_3=-c_3,\quad \mu_3=b_3-b_6,\\
	&\lambda_4=-c_4,\quad \mu_4=b_4,\quad \lambda_5=-c_5,\quad \mu_5=b_5-b_7-1.
\end{split}\]
Those variables satisfy a rational Hamiltonian system of fifth order with a symplectic form $\omega=\sum_{j=1}^{5}d\lambda_j\wedge d\mu_j$; we do not give its explicit formula here.
Then the third and fourth relation of \eqref{Eq:Residue_42_33_33_222_7th} are described as
\[
	\lambda_1\mu_1 + \mu_1\mu_5 - \theta_2 - \rho_2 = 0,\quad
	\lambda_1\mu_1 - \lambda_3\mu_3 - \lambda_5\mu_5 - \theta_1 - \theta_2 - 2\rho_2 = 0.
\]
Substituting them to the symplectic form $\omega$ again, we obtain
\[\begin{split}
	\omega &= d\lambda_1\wedge d\mu_1 + d\lambda_2\wedge d\mu_2 + d\lambda_3\wedge d\frac{1}{\lambda_3}\left(\lambda_1\mu_1+\lambda_1\lambda_5-\frac{(\theta_2+\rho_2)\lambda_5}{\mu_1}-\theta_1-\theta_2-2\rho_2\right)\\
	&\quad + d\lambda_4\wedge d\mu_4 + d\lambda_5\wedge d\left(-\lambda_1+\frac{\theta_2+\rho_2}{\mu_1}\right)\\
	&= d\lambda_3\left(\lambda_1-\frac{\theta_2+\rho_2}{\mu_1}\right)\wedge d\frac{\mu_1+\lambda_5}{\lambda_3} + d\lambda_2\wedge d\mu_2 + d\lambda_4\wedge d\mu_4.
\end{split}\]
Hence we can take
\[
	q_1 = \frac{\mu_1+\lambda_5}{\lambda_3},\quad
	p_1 = -\lambda_3\left(\lambda_1-\frac{\theta_2+\rho_2}{\mu_1}\right),\quad
	q_2 = \lambda_4,\quad
	p_2 = \mu_4,\quad
	\frac{q_3}{t} = -\lambda_2,\quad
	tp_3 = -\mu_2.
\]
Then, in a similar manner as Section \ref{Sec:Derivation_21_21_21_21_111}, we arrive at
\begin{thm}
The dependent variables $p_j,q_j$ $(j=1,2,3)$ satisfy the system $\mathcal{H}^{42,33,33,222}$.
\end{thm}

In the above, we used 31 relations of \eqref{Eq:Residue_42_33_33_222} to derive the system $\mathcal{H}^{42,33,33,222}$.
And the rest 4 relations have not been used yet.
They are used to determine the variables $b_6,b_7,c_3,c_5$ as rational expressions in $(q_j,p_j)$.
Hence we can show that the components of the matrices $\widetilde{A}_1,\widetilde{A}_2,\widetilde{A}_3$ are rational in $(q_j,p_j)$; we do not give their explicit formulas here.

\subsection{Spectral type $51,33,222,222$}

In this case, we consider a Fuchsian system
\[
	\frac{\partial}{\partial x}Y(x) = \left(\frac{A_1}{x-t}+\frac{A_2}{x-1}+\frac{A_3}{x}\right)Y(x),
\]
with a Riemann scheme
\[
	\left\{\begin{array}{cccc}
		x=t & x=1 & x=0 & x=\infty \\
		\theta_1 & \theta_2 & \theta_{3,1} & \rho_1 \\
		\theta_1 & 0 & \theta_{3,1} & \rho_1 \\
		\theta_1 & 0 & \theta_{3,2} & \rho_2 \\
		0 & 0 & \theta_{3,2} & \rho_2 \\
		0 & 0 & 0 & \rho_3 \\
		0 & 0 & 0 & \rho_3 \\
	\end{array}\right\}.
\]
Note that a Fuchsian relation $3\theta_1+\theta_2+2\theta_{3,1}+2\theta_{3,2}+2\rho_1+2\rho_2+2\rho_3=0$ is satisfied.

In a similar manner as Section \ref{Sec:Derivation_21_21_21_21_111}, the residue matrices are transformed to
\[\begin{split}
	&\widetilde{A}_1 = \begin{pmatrix}I_3\\\widetilde{B}_1\end{pmatrix}\begin{pmatrix}\theta_1I_3-\widetilde{C}_1\widetilde{B}_1&\widetilde{C}_1\end{pmatrix},\quad
	\widetilde{B}_1 = \begin{pmatrix}a_1&a_2&0\\a_3&b_1&b_2\\a_4&a_5&b_3\end{pmatrix},\quad
	\widetilde{C}_1 = \begin{pmatrix}0&0&1\\1&c_1&0\\a_6&c_2&c_3\end{pmatrix},\\
	&\widetilde{A}_2 = \begin{pmatrix}1\\a_7\\a_8\\a_9\\a_{10}\\a_{11}\end{pmatrix}\begin{pmatrix}\theta_2-a_7-a_8-a_{10}&1&1&0&1&0\end{pmatrix},\quad
	\widetilde{A}_3 = \begin{pmatrix}\theta_{3,1}&0&a_{12}&a_{13}&a_{14}&a_{15}\\0&\theta_{3,1}&a_{16}&a_{17}&a_{18}&a_{19}\\0&0&\theta_{3,2}&0&a_{20}&a_{21}\\0&0&0&\theta_{3,2}&a_{22}&a_{23}\\0&0&0&0&0&0\\0&0&0&0&0&0\end{pmatrix},
\end{split}\]
and
\begin{equation}\label{Eq:Residue_51_33_222_222}
	\widetilde{A}_{\infty} := -\sum_{i=1}^{3}\widetilde{A}_i = \begin{pmatrix}\rho_1&0&0&0&0&0\\0&\rho_1&0&0&0&0\\a_{24}&a_{25}&\rho_2&0&0&0\\a_{26}&a_{27}&0&\rho_2&0&0\\a_{28}&a_{29}&a_{30}&a_{31}&\rho_3&0\\a_{32}&a_{33}&a_{34}&a_{35}&0&\rho_3\end{pmatrix}.
\end{equation}
By using \eqref{Eq:Residue_51_33_222_222}, we can show that the variables $a_1,\ldots.a_{35}$ are rational in $(b_j,c_j)$; we do not give their explicit formulas here.
Then the dependent variables $b_j,c_j$ $(j=1,2,3)$ satisfy a Hamiltonian system
\begin{equation}\label{Eq:Sch_Ham_51_33_222_222}
	\frac{\partial b_j}{\partial t} = \{H,b_j\},\quad
	\frac{\partial c_j}{\partial t} = \{H,c_j\},\quad
	H = \frac{\mathrm{tr}\widetilde{A}_1\widetilde{A}_2}{t-1} + \frac{\mathrm{tr}\widetilde{A}_1\widetilde{A}_3}{t},
\end{equation}
with a symplectic form $\omega=\sum_{j=1}^{3}db_j\wedge dc_j$.
Note that the Hamiltonian $H$ turns out to be a polynomial in $(b_j,c_j)$, although the components of the matrices $\widetilde{A}_1,\widetilde{A}_2,\widetilde{A}_3$ are rational.

Under the system \eqref{Eq:Sch_Ham_51_33_222_222}, we consider a canonical transformation
\[\begin{split}
	&q_1 = -b_2,\quad
	p_1 = -c_2,\quad
	q_2 = \frac{1}{b_3},\quad
	p_2 = -b_3\left(b_3c_3-\frac{\theta_1}{2}-\frac{\theta_2}{2}-\theta_{3,2}-\rho_2\right),\\
	&\frac{q_3}{t} = -\frac{1}{b_1},\quad
	tp_3 = b_1\left(b_1c_1+\frac{\theta_1}{2}+\frac{\theta_2}{2}+\rho_3\right).
\end{split}\]
Then, by a direct computation, we arrive at
\begin{thm}
The dependent variables $p_j,q_j$ $(j=1,2,3)$ satisfy the system $\mathcal{H}^{51,33,222,222}$.
\end{thm}

\section{Particular Solutions}\label{Sec:Solutions}

In this section, we give particular solutions of the six-dimensional Painlev\'{e} systems in terms of the hypergeometric functions.

\subsection{Spectral type $21,21,21,21,111$}

Under the system $\mathcal{H}^{21,21,21,21,111}$, we consider a specialization
\[
	p_1 - \frac{\alpha_1}{q_1} = p_2 = q_3 = \alpha_0 + \alpha_1 + \alpha_5 + 1 = 0.
\]
Also we set
\[
	\frac{y_1}{y_0} = q_1,\quad
	\frac{y_2}{y_0} = q_2,\quad
	\frac{y_3}{y_0} = -t_2q_1p_3,
\]
where the variable $y_0$ satisfies a Pfaff system
\[
	t_1(t_1-1)\frac{\partial}{\partial t_1}\log y_0 = (\alpha_5+1)q_1 + \alpha_1t_1,\quad
	t_2(t_2-1)\frac{\partial}{\partial t_2}\log y_0 = t_2q_1p_3 + (\alpha_5-\rho_3+1)q_2 - \alpha_2t_2.
\]
Then we have
\begin{thm}
A vector of variables $\mathbf{y}={}^t[y_0,y_1,y_2,y_3]$ satisfies a rigid system
\begin{equation}\label{Eq:HGE_31_31_22_211}
	\frac{\partial\mathbf{y}}{\partial t_i} = \left(\frac{M^{(i)}_t}{t_i-t_{i+1}}+\frac{M^{(i)}_1}{t_i-1}+\frac{M^{(i)}_0}{t_i}\right)\mathbf{y}\quad (i\in\mathbb{Z}/2\mathbb{Z}),
\end{equation}
with matrices
\[\begin{split}
	&M^{(1)}_t = \begin{pmatrix}0&0&0&0\\0&-\alpha_2&-\alpha_1&-1\\0&\alpha_2&\alpha_1&1\\0&-\alpha_2\rho_3&-\alpha_1\rho_3&-\rho_3\end{pmatrix},\quad
	M^{(1)}_1 = \begin{pmatrix}\alpha_1&\alpha_5+1&0&0\\\alpha_1&\alpha_5+1&0&0\\0&0&0&0\\0&0&\alpha_1\rho_3&\alpha_1+\alpha_5+1\end{pmatrix},\\
	&M^{(1)}_0 = \begin{pmatrix}0&-\alpha_5-1&0&0\\0&\alpha_1-\alpha_3+1&0&0\\0&-\alpha_2&0&0\\0&\alpha_2\rho_3&0&0\end{pmatrix},
\end{split}\]
and
\[\begin{split}
	&M^{(2)}_t = \begin{pmatrix}0&0&0&0\\0&-\alpha_2&-\alpha_1&-1\\0&\alpha_2&\alpha_1&1\\0&-\alpha_2\rho_3&-\alpha_1\rho_3&-\rho_3\end{pmatrix},\quad
	M^{(2)}_1 = \begin{pmatrix}-\alpha_2&0&\alpha_5-\rho_3+1&-1\\0&0&0&0\\-\alpha_2&0&\alpha_5-\rho_3+1&-1\\0&0&0&0\end{pmatrix},\\
	&M^{(2)}_0 = \begin{pmatrix}0&0&-\alpha_5+\rho_3-1&1\\0&0&\alpha_1&1\\0&0&-\alpha_2-\alpha_3+\rho_3+1&0\\0&0&0&-\alpha_2-\alpha_3+\rho_3+1\end{pmatrix}.
\end{split}\]
\end{thm}

The Riemann scheme of the system \eqref{Eq:HGE_31_31_22_211} is given by
\[
	\left\{\begin{array}{cccc}
		t_1=t_2 & t_1=1 & t_1=0 & t_1=\infty \\
		\alpha_1-\alpha_2-\rho_3 & \alpha_1+\alpha_5+1 & \alpha_1-\alpha_3+1 & -\alpha_1+\alpha_2+\alpha_3-\alpha_5-2 \\
		0 & \alpha_1+\alpha_5+1 & 0 & -\alpha_1-\alpha_5+\rho_3-1 \\
		0 & 0 & 0 & -\alpha_1 \\
		0 & 0 & 0 & -\alpha_1
	\end{array}\right\},
\]
and
\[
	\left\{\begin{array}{cccc}
		t_2=t_1 & t_2=1 & t_2=0 & t_2=\infty \\
		\alpha_1-\alpha_2-\rho_3 & -\alpha_2+\alpha_5-\rho_3+1 & -\alpha_2-\alpha_3+\rho_3+1 & -\alpha_1+\alpha_2+\alpha_3-\alpha_5-2 \\
		0 & 0 & -\alpha_2-\alpha_3+\rho_3+1 & \alpha_2+\alpha_3-1 \\
		0 & 0 & 0 & \alpha_2 \\
		0 & 0 & 0 & \alpha_2
	\end{array}\right\}.
\]
Namely, for each $i\in\mathbb{Z}/2\mathbb{Z}$, the system \eqref{Eq:HGE_31_31_22_211} is a Fuchsian one with a spectral type $31,31,22,211$.

\subsection{Spectral type $31,31,22,22,22$}

Under the system $\mathcal{H}^{31,31,22,22,22}$, we consider a specialization
\[
	q_1 = q_2 = q_3 = \alpha_1 = 0.
\]
Note that such a specialization implies $\widetilde{A}_1=O$.
We also set
\[
	\frac{y_1}{y_0} = -t_1p_1,\quad
	\frac{y_2}{y_0} = -t_2p_2,\quad
	\frac{y_3}{y_0} = t_1t_2(p_1p_2-\alpha_2p_3),
\]
where the variable $y_0$ satisfies a Pfaff system
\[
	t_i(t_i-1)\frac{\partial}{\partial t_i}\log y_0 = t_ip_i - (\alpha_0+\alpha_5+1)t_i - \alpha_3(t_i-1)\quad (i\in\mathbb{Z}/2\mathbb{Z}).
\]
Then we have
\begin{thm}
A vector of variables $\mathbf{y}={}^t[y_0,y_1,y_2,y_3]$ satisfies a rigid system
\begin{equation}\label{Eq:HGE_31_22_22_22}
	\frac{\partial\mathbf{y}}{\partial t_i} = \left(\frac{M^{(i)}_t}{t_i-t_{i+1}}+\frac{M^{(i)}_1}{t_i-1}+\frac{M^{(i)}_0}{t_i}\right)\mathbf{y}\quad (i\in\mathbb{Z}/2\mathbb{Z}),
\end{equation}
with matrices
\[\begin{split}
	&M^{(1)}_t = \begin{pmatrix}0&0&0&0\\0&\alpha_2&-\alpha_2&0\\0&-\alpha_2&\alpha_2&0\\0&0&0&0\end{pmatrix},\\
	&M^{(1)}_1 = \begin{pmatrix}-(\alpha_0+\alpha_5+1)&-1&0&0\\\alpha_0(\alpha_0+\alpha_5+1)&\alpha_0&0&0\\0&0&-(\alpha_0+\alpha_2+\alpha_5+1)&-1\\0&0&(\alpha_0+\alpha_2)(\alpha_0+\alpha_2+\alpha_5+1)&\alpha_0+\alpha_2\end{pmatrix},\\
	&M^{(1)}_0 = \begin{pmatrix}-\alpha_3&1&0&0\\0&0&0&0\\0&\alpha_2&-\alpha_3&1\\0&0&0&0\end{pmatrix},
\end{split}\]
and
\[
	M^{(2)}_t = E_{23}M^{(1)}_tE_{23},\quad
	M^{(2)}_1 = E_{23}M^{(1)}_1E_{23},\quad
	M^{(2)}_0 = E_{23}M^{(1)}_0E_{23},\quad
	E_{23} = \begin{pmatrix}1&0&0&0\\0&0&1&0\\0&1&0&0\\0&0&0&1\end{pmatrix}.
\]
\end{thm}

The Riemann scheme of the system \eqref{Eq:HGE_31_22_22_22} is given by
\[
	\left\{\begin{array}{cccc}
		t_1=t_2 & t_1=1 & t_1=0 & t_1=\infty \\
		2\alpha_2 & -\alpha_5-1 & -\alpha_3 & \alpha_0+\alpha_3+\alpha_5+1 \\
		0 & -\alpha_5-1 & -\alpha_3 & \alpha_0+\alpha_3+\alpha_5+1 \\
		0 & 0 & 0 & -\alpha_0-\alpha_2 \\
		0 & 0 & 0 & -\alpha_0-\alpha_2
	\end{array}\right\},
\]
and
\[
	\left\{\begin{array}{cccc}
		t_2=t_1 & t_2=1 & t_2=0 & t_2=\infty \\
		2\alpha_2 & -\alpha_5-1 & -\alpha_3 & \alpha_0+\alpha_3+\alpha_5+1 \\
		0 & -\alpha_5-1 & -\alpha_3 & \alpha_0+\alpha_3+\alpha_5+1 \\
		0 & 0 & 0 & -\alpha_0-\alpha_2 \\
		0 & 0 & 0 & -\alpha_0-\alpha_2
	\end{array}\right\}.
\]
Namely, for each $i\in\mathbb{Z}/2\mathbb{Z}$, the system \eqref{Eq:HGE_31_22_22_22} is a Fuchsian one with a spectral type $31,22,22,22$.

\subsection{Spectral type $21,111,111,111$}

Under the system $\mathcal{H}^{21,111,111,111}$, we consider a specialization
\[
	p_1 = p_2 = p_3 = \eta = 0.
\]
Also we set
\[
	\frac{y_1}{y_0} = \frac{q_1}{t},\quad
	\frac{y_2}{y_0} = \frac{q_2}{t},\quad
	\frac{y_3}{y_0} = q_3,
\]
where the variable $y_0$ satisfies a Pfaff system
\[
	t(t-1)\frac{d}{dt}\log y_0 = -\alpha_1q_1 - \alpha_5q_2 -\alpha_3t.
\]
Then we have
\begin{thm}
A vector of variables $\mathbf{y}={}^t[y_0,y_1,y_2,y_3]$ satisfies a rigid system
\begin{equation}\label{Eq:HGE_211_211_211}
	\frac{d\mathbf{y}}{dt} = \left(\frac{M_1}{t-1}+\frac{M_0}{t}\right)\mathbf{y},
\end{equation}
with matrices
\[
	M_1 = \begin{pmatrix}-\alpha_3&-\alpha_1&-\alpha_5&0\\-\alpha_3&-\alpha_1+\theta_{2,1}&-\alpha_5-\theta_{2,1}&\alpha_3\\-\alpha_3&-\alpha_1&-\alpha_5&0\\0&-\theta_{2,1}&\theta_{2,1}&-\alpha_3\end{pmatrix},\quad
	M_0 = \begin{pmatrix}0&0&0&0\\\alpha_3&-\alpha_2-\alpha_3&0&-\alpha_3\\\alpha_3&\alpha_1&\alpha_4+\alpha_5-1&0\\0&0&0&0\end{pmatrix}.
\]
\end{thm}

The Riemann scheme of the system \eqref{Eq:HGE_211_211_211} is given by
\[
	\left\{\begin{array}{ccc}
		t=1 & t=0 & t=\infty \\
		-\alpha_1-\alpha_3-\alpha_5 & \alpha_4+\alpha_5-1 & -\alpha_4+1 \\
		-\alpha_3+\theta_{2,1} & -\alpha_2-\alpha_3 & \alpha_1+\alpha_2+\alpha_3-\theta_{2,1} \\
		0 & 0 & \alpha_3 \\
		0 & 0 & \alpha_3
	\end{array}\right\}.
\]
Namely, the system \eqref{Eq:HGE_211_211_211} is a Fuchsian one with a spectral type $211,211,211$.

\subsection{Spectral type $31,22,211,1111$}

Under the system $\mathcal{H}^{31,22,211,1111}$, we consider a specialization
\[
	q_1p_1-\alpha_1 = q_2 = q_3 = \alpha_1+\alpha_3-\eta = 0.
\]
Note that such a specialization implies $\widetilde{A}_1=O$.
We also set
\[
	\frac{y_1}{y_0} = tp_1,\quad
	\frac{y_2}{y_0} = tp_2 - \frac{t^2p_1p_3}{\eta},\quad
	\frac{y_3}{y_0} = \frac{t^2p_1p_3}{\eta}.
\]
where the variable $y_0$ satisfies a Pfaff system
\[
	t(t-1)\frac{d}{dt}\log y_0 = tp_1 + tp_2 + (\alpha_1-\alpha_5-\rho_4)t - (\alpha_0+\alpha_1+\alpha_2)(t-1)
\]
Then we have
\begin{thm}
A vector of variables $\mathbf{y}={}^t[y_0,y_1,y_2,y_3]$ satisfies a rigid system
\begin{equation}\label{Eq:HGE_22_211_1111}
	\frac{d\mathbf{y}}{dt} = \left(\frac{M_1}{t-1}+\frac{M_0}{t}\right)\mathbf{y},
\end{equation}
with matrices
\[\begin{split}
	&M_1 = \begin{pmatrix}\alpha_1-\alpha_5-\rho_4&-1&-1&-1\\\alpha_1(\alpha_1-\eta)&-\alpha_1-\alpha_5+\eta-\rho_4&-\alpha_1&-\alpha_1+\eta\\-\frac{(\alpha_1-\eta)(\alpha_5+\rho_4)(\eta-\rho_4)}{\eta}&\frac{(\alpha_5+\rho_4)(\eta-\rho_4)}{\eta}&\eta-\rho_4&0\\-\frac{\alpha_1(\alpha_5-\eta+\rho_4)\rho_4}{\eta}&\frac{(\alpha_5-\eta+\rho_4)\rho_4}{\eta}&0&-\rho_4\end{pmatrix},\\
	&M_0 = \begin{pmatrix}-\alpha_0-\alpha_1-\alpha_2&1&1&1\\0&-\alpha_0&\alpha_1&\alpha_1-\eta\\0&0&0&0\\0&0&0&0\end{pmatrix}.
\end{split}\]
\end{thm}

The Riemann scheme of the system \eqref{Eq:HGE_22_211_1111} is given by
\[
	\left\{\begin{array}{ccc}
		t=1 & t=0 & t=\infty \\
		-\alpha_5+\eta-2\rho_4 & -\alpha_0-\alpha_1-\alpha_2 & \alpha_0+\alpha_2+\alpha_5+\rho_4 \\
		-\alpha_5+\eta-2\rho_4 & -\alpha_0 & \alpha_0+\alpha_1+\alpha_5-\eta+\rho_4 \\
		0 & 0 & -\eta+\rho_4 \\
		0 & 0 & \rho_4
	\end{array}\right\}.
\]
Namely, the system \eqref{Eq:HGE_22_211_1111} is a Fuchsian one with a spectral type $22,211,1111$.

\subsection{Remark: Theorycal Background}

Let $\mathbf{m}=\{m_2,\ldots,m_{N+3}\}$, where $m_i=(m_{i,1},\ldots,m_{i,l_i})$, be a $(N+2)$-tuples of partitions of natural number $L$ such that
\[
	NL^2 - \sum_{i=2}^{N+3}\sum_{j=1}^{l_i}m_{i,j}^2 + 2 = 0.
\]
Also let $\widetilde{\mathbf{m}}=\{m_1,m_2,\ldots,m_{N+3}\}$, where $m_1=(m_{1,1},m_{1,2})$, be a $(N+3)$-tuples of partitions of $L$.
Note that a Fuchsian system with a spectral type $\mathbf{m}$ (or resp. $\widetilde{\mathbf{m}}$) contains 0 (or resp. $2m_{1,1}m_{1,2}$) accessory parameters.
We consider a Schlesinger system \eqref{Eq:Sch_Ham} associated with a spectral type $\widetilde{\mathbf{m}}$, which is rewritten to
\begin{equation}\begin{split}\label{Eq:Sch_Ham_HGE}
	\frac{\partial B_j}{\partial t_i} = \frac{A_i}{t_i-t_j}B_j,\quad
	\frac{\partial C_j}{\partial t_i} = -C_j\frac{A_i}{t_i-t_j},\quad
	\frac{\partial B_i}{\partial t_i} = -\sum_{j=1; j\neq i}^{N+2}\frac{A_j}{t_i-t_j}B_i,\quad
	\frac{\partial C_i}{\partial t_i} = C_i\sum_{j=1; j\neq i}^{N+2}\frac{A_j}{t_i-t_j}\\
	(i=1,\ldots,N; j=1,\ldots,N+2; j\neq i).
\end{split}\end{equation}
Such a system admits a particular solution given by a rigid system.
\begin{lem}
The system \eqref{Eq:Sch_Ham_HGE} admits a specialization $C_1=O$.
Then a matrix of variables $B_1$ satisfies
\[
	\frac{\partial B_1}{\partial t_1} = -\sum_{j=2}^{N+2}\frac{A_j}{t_1-t_j}B_1,\quad
	\frac{\partial B_1}{\partial t_i} = \frac{A_j}{t_i-t_1}B_1\quad (i=2,\ldots,N).
\]
\end{lem}

In the previous subsections, we have given specializations and Pfaff systems in order to derive rigid systems.
Their origins can be clarified with the aid of this lemma.
Furthermore, this lemma suggests that we can always give a particular solutions of a Painlev\'{e} system $\mathcal{H}^{\widetilde{\mathbf{m}}}$ by a rigid system with a spectral type $\mathbf{m}$.

\appendix

\section{Hamiltonians of the six-dimensional Painlev\'{e} system}\label{App:Hamiltonian_6dim}

In this section, we recall the Hamiltonian of the six-dimensional Painlev\'{e} systems which have been already derived in \cite{G,T1,Kaw,FISS}.

\subsection{Spectral type $11,11,11,11,11,11$}

\[\begin{split}
	H^{11,11,11,11,11,11}_i &= H_{\rm{VI}}(\alpha_6+1,\alpha_{i+1}+\alpha_{i-1}+\alpha_4,\alpha_5,\alpha_i;q_i,p_i;t_i)\\
	&\quad + q_iq_{i+1}p_{i+1}(2q_ip_i+q_{i+1}p_{i+1}+2\alpha_0+\alpha_6+1)\\
	&\quad - \frac{1}{t_i-t_{i+1}}q_iq_{i+1}\{t_i(t_i-1)p_i^2-2t_i(t_{i+1}-1)p_ip_{i+1}+(t_i-1)t_{i+1}p_{i+1}^2\}\\
	&\quad + \alpha_i\frac{t_i}{t_i-t_{i+1}}q_{i+1}\{(t_i-1)p_i-(t_{i+1}-1)p_{i+1}\} - \alpha_{i+1}\frac{(t_i-1)t_{i+1}}{t_i-t_{i+1}}q_i(p_i-p_{i+1})\\
	&\quad + q_iq_{i-1}p_{i-1}(2q_ip_i+q_{i-1}p_{i-1}+2\alpha_0+\alpha_6+1)\\
	&\quad - \frac{1}{t_i-t_{i-1}}q_iq_{i-1}\{t_i(t_i-1)p_i^2-2t_i(t_{i-1}-1)p_ip_{i-1}+(t_i-1)t_{i-1}p_{i-1}^2\}\\
	&\quad + \alpha_i\frac{t_i}{t_i-t_{i-1}}q_{i-1}\{(t_i-1)p_i-(t_{i-1}-1)p_{i-1}\} - \alpha_{i-1}\frac{(t_i-1)t_{i-1}}{t_i-t_{i-1}}q_i(p_i-p_{i-1}),
\end{split}\]
for $i\in\mathbb{Z}/3\mathbb{Z}$, where $2\alpha_0+\alpha_1+\alpha_2+\alpha_3+\alpha_4+\alpha_5+\alpha_6=0$.

\subsection{Spectral type $31,31,1111,1111$}

\[\begin{split}
	H^{31,31,1111,1111} &= H_{\rm{VI}}(-\alpha_1+\eta,\alpha_2+\alpha_4+\alpha_6,\alpha_0,\alpha_3+\alpha_5+\alpha_7-\eta;q_i,p_i;t)\\
	&\quad + H_{\rm{VI}}(-\alpha_3+\eta,\alpha_4+\alpha_6,\alpha_0+\alpha_2,\alpha_1+\alpha_5+\alpha_7-\eta;q_i,p_i;t)\\
	&\quad + H_{\rm{VI}}(-\alpha_5+\eta,\alpha_6,\alpha_0+\alpha_2+\alpha_4,\alpha_1+\alpha_3+\alpha_7-\eta;q_i,p_i;t)\\
	&\quad + (q_1-1)(q_2-t)\{(q_1p_1+\alpha_1)p_2+p_1(p_2q_2+\alpha_3)\}\\
	&\quad + (q_1-1)(q_3-t)\{(q_1p_1+\alpha_1)p_3+p_1(p_3q_3+\alpha_5)\}\\
	&\quad + (q_2-1)(q_3-t)\{(q_2p_2+\alpha_3)p_3+p_2(p_3q_3+\alpha_5)\},
\end{split}\]
where $\alpha_0+\alpha_1+\alpha_2+\alpha_3+\alpha_4+\alpha_5+\alpha_6+\alpha_7=1$.

\subsection{Spectral type $33,33,33,321$}

\[\begin{split}
	H^{33,33,33,321} &= \mathrm{tr}\{Q(Q-1)P(Q-t)P-(\alpha_1-1)Q(Q-1)P\\
	&\quad -\alpha_3Q(Q-t)P-\alpha_4(Q-1)(Q-t)P+\alpha_2(\alpha_0+\alpha_2)Q\}
\end{split}\]
where $\alpha_0+\alpha_1+2\alpha_2+\alpha_3+\alpha_4=1$ and
\[\begin{split}
	P &= -\frac{1}{t}\begin{pmatrix}\frac{1}{3}p_1+q_1^2p_3&\frac{2}{3}q_2p_2+(q_3-q_1^3)p_3+2\alpha_2+2\alpha_5&\frac{1}{3}(q_3-q_1^3)p_2+2q_2^2p_3\\\frac{1}{3}p_2&\frac{1}{3}p_1+(q_2+q_1^2)p_3&\frac{1}{3}q_2p_2+(q_3-q_1^3)p_3+\alpha_2+\alpha_5\\p_3&\frac{1}{3}p_2&\frac{1}{3}p_1-(q_2-q_1^2)p_3\end{pmatrix},\\
	Q &= -t\begin{pmatrix}q_1&2q_2&q_3-q_1^3\\1&q_1&q_2\\0&1&q_1\end{pmatrix}.
\end{split}\]

\subsection{Spectral type $51,33,33,111111$}

\[\begin{split}
	H^{51,33,33,111111} &= H_{\rm{VI}}(\alpha_0,\alpha_1,\alpha_3+2\alpha_4+\alpha_5+2\alpha_6+\alpha_7,\alpha_3+\alpha_5+\alpha_8;q_1,p_1;t)\\
	&\quad + H_{\rm{VI}}(\alpha_0+2\alpha_2+\alpha_3,\alpha_1+\alpha_3,2\alpha_4+\alpha_5+2\alpha_6+\alpha_7,\alpha_5+\alpha_8;q_2,p_2;t)\\
	&\quad + H_{\rm{VI}}(\alpha_0+2\alpha_2+\alpha_3+2\alpha_4+\alpha_5,\alpha_1+\alpha_3+\alpha_5,2\alpha_6+\alpha_7,\alpha_8;q_3,p_3;t)\\
	&\quad + 2(q_1-t)p_1q_2\{(q_2-1)p_2+\alpha_4\} + 2(q_1-t)p_1q_3\{(q_3-1)p_3+\alpha_6\}\\
	&\quad + 2(q_2-t)p_2q_3\{(q_3-1)p_3+\alpha_6\},
\end{split}\]
where $\alpha_0+\alpha_1+2\alpha_2+2\alpha_3+2\alpha_4+2\alpha_5+2\alpha_6+\alpha_7+\alpha_8=1$.

\section{Hamiltonians of the four-dimensional Painlev\'{e} system}\label{App:Hamiltonian_4dim}

In this section, we recall the Hamiltonian of the four-dimensional Painlev\'{e} systems which have been classified by Sakai \cite{Sak}.

\subsection{Spectral type $11,11,11,11,11$}

\[\begin{split}
	H^{11,11,11,11,11}_i &= H_{\rm{VI}}(\alpha_5+1,\alpha_{i+1}+\alpha_3,\alpha_4,\alpha_i;q_i,p_i;t_i) + q_1q_2p_{i+1}(2q_ip_i+q_{i+1}p_{i+1}+2\alpha_0+\alpha_5+1)\\
	&\quad - \frac{1}{t_i-t_{i+1}}q_1q_2\{t_i(t_i-1)p_i^2-2t_i(t_{i+1}-1)p_1p_2+(t_i-1)t_{i+1}p_{i+1}^2\}\\
	&\quad + \alpha_i\frac{t_i}{t_i-t_{i+1}}q_{i+1}\{(t_i-1)p_i-(t_{i+1}-1)p_{i+1}\} - \alpha_{i+1}\frac{(t_i-1)t_{i+1}}{t_i-t_{i+1}}q_i(p_i-p_{i+1}),
\end{split}\]
for $i\in\mathbb{Z}/2\mathbb{Z}$, where $2\alpha_0+\alpha_1+\alpha_2+\alpha_3+\alpha_4+\alpha_5=0$.

\subsection{Spectral type $21,21,111,111$}

\[\begin{split}
	H^{21,21,111,111} &= H_{\rm{VI}}(-\alpha_1+\eta,\alpha_2,\alpha_0+\alpha_4,\alpha_3+\alpha_5-\eta;q_1,p_1;t)\\
	&\quad + H_{\rm{VI}}(-\alpha_5+\eta,\alpha_0+\alpha_2,\alpha_4,\alpha_1+\alpha_3-\eta;q_2,p_2;t)\\
	&\quad + (q_1-t)(q_2-1)\{(q_1p_1+\alpha_1)p_2+p_1(q_2p_2+\alpha_5)\},
\end{split}\]
where $\alpha_0+\alpha_1+\alpha_2+\alpha_3+\alpha_4+\alpha_5=1$.

\subsection{Spectral type $22,22,22,211$}

\[\begin{split}
	H^{22,22,22,211} &= \mathrm{tr}\{Q(Q-1)P(Q-t)P-(\alpha_1-1)Q(Q-1)P\\
	&\quad -\alpha_3Q(Q-t)P-\alpha_4(Q-1)(Q-t)P+\alpha_2(\alpha_0+\alpha_2)Q\},
\end{split}\]
where $\alpha_0+\alpha_1+2\alpha_2+\alpha_3+\alpha_4=1$ and
\[
	P = \frac{1}{t}\begin{pmatrix}\frac{1}{2}p_1&-p_2\\q_2p_2+\alpha_2+\alpha_5&\frac{1}{2}p_1\end{pmatrix},\quad
	Q = t\begin{pmatrix}-q_1&-1\\q_2&-q_1\end{pmatrix}.
\]

\subsection{Spectral type $31,22,22,1111$}

\[\begin{split}
	H^{31,22,22,1111} &= H_{\rm{VI}}(\alpha_0,\alpha_1,\alpha_3+2\alpha_4+\alpha_5,\alpha_3+\alpha_6;q_1,p_1;t)\\
	&\quad + H_{\rm{VI}}(\alpha_0+2\alpha_2+\alpha_3,\alpha_1+\alpha_3,\alpha_5,\alpha_6;q_2,p_2;t)\\
	&\quad + 2(q_1-t)p_1q_2\{(q_2-1)p_2+\alpha_4\},
\end{split}\]
where $\alpha_0+\alpha_1+2\alpha_2+2\alpha_3+2\alpha_4+\alpha_5+\alpha_6=1$.

\section*{Acknowledgement}
The author would like to express his gratitude to the collaborator in the previous work \cite{FISS}, Dr. Kenta Fuji, Mr. Keisuke Inoue and Mr. Keisuke Shinomiya.
The auther is also grateful to Professors Kazuki Hiroe, Hiroshi Kawakami, Hajime Nagoya, Masatoshi Noumi, Toshio Oshima, Hidetaka Sakai, Teruhisa Tsuda and Yasuhiko Yamada for valuable discussions and advices.


\end{document}